\documentclass[12pt,draft,leqno]{article}
\usepackage{amssymb, eucal, latexsym}

0in \oddsidemargin 1,2cm \evensidemargin 0in

\newtheorem{theorem}{Theorem}[section]
\newtheorem{lm}[theorem]{Lemma}
\newtheorem{exa}[theorem]{Example}

\newtheorem{cor}[theorem]{Corollary}
\newtheorem{pro}[theorem]{Proposition}
\newtheorem{defi}[theorem]{Definition}

\newtheorem{nota}[theorem]{Notation}
\newtheorem{notas}[theorem]{Notations}
\newtheorem{rem}[theorem]{Remark}

\newtheorem{fact}[theorem]{Fact}

\def\p{\varphi}
\def\a{\alpha}

\def\d{\delta}
\def\ep{\varepsilon}
\def\g{\gamma}
\def\GA{\Gamma}

\def\l{\lambda}
\def\LAM{\Lambda}

\def\TE{\Theta}

\def\s{\sigma}
\def\SI{\Sigma}

\def\di{\diamond}

\def\fs{\hat{f}}
\def\ps{\hat{\varphi}}
\def\lag{\lambda_A^g}
\def\lbg{\lambda_B^g}

\def\lbc{\lambda_B^C}

\def\lra{\longrightarrow}

\def\dar{\downarrow}

\def\sbe{\subseteq}
\def\spe{\supseteq}
\def\stm{\setminus}
\def\ems{\emptyset}
\def\nes{\neq\emptyset}

\def\fa{\forall}

\def\we{\wedge}

\def\bv{\bigvee}

\def\ap{^\prime}
\def\inv{^{-1}}
\def\st{\ |\ }

\def\llx{\ll_{\rho}}
\def\lle{\ll_{\eta}}

\def\nin{\not\in}

\def\card #1{\vert #1 \vert}

\def\AA{{\cal A}}
\def\BB{{\cal B}}
\def\CC{{\cal C}}
\def\DD{{\cal D}}
\def\EE{{\cal E}}

\def\KK{{\cal K}}

\def\MM{{\cal M}}
\def\OO{{\cal O}}
\def\PP{{\cal P}}
\def\QQ{{\cal Q}}
\def\RR{{\cal R}}
\def\TT{{\cal T}}
\def\SS{{\cal S}}

\def\ZZ{{\cal Z}}

\def\Bo{{\bf Bool}}

\def\HLC{{\bf HLC}}
\def\GBPL{{\bf GBPL}}
\def\LBA{{\bf LBA}}
\def\DLC{{\bf DLC}}

\def\DZLC{{\bf DZLC}}
\def\DPZLC{{\bf DPZLC}}

\def\PAL{{\bf PAL}}

\def\ZHC{{\bf ZHC}}
\def\PZLC{{\bf PZLC}}
\def\ZLC{{\bf ZLC}}
\def\PLBA{{\bf PLBA}}
\def\PZLBA{{\bf PZLBA}}
\def\ZLBA{{\bf ZLBA}}

\def\Bool{{\bf Bool}}

\def\B{\mbox{{\boldmath $B$}}}

\def\2{\mbox{{\bf 2}}}
\def\3{\mbox{{\bf 3}}}

\def\int{\mbox{{\rm int}}}

\def\cl{\mbox{{\rm cl}}}

\def\doc{\hspace{-1cm}{\em Proof.}~~}
\def\sq{\hspace*{\fill} \hbox{\vrule\vbox{\hrule\phantom{o}\hrule}\vrule}}
\def\sqs{\sq \vspace{2mm}}

\def\scx{u_x^C}

\def\scy{u_y^C}

\def\tcx{t_X^C}
\def\tcy{t_Y^C}

\def\BBBB{{\rm I}\!{\rm B}}

\title{{\LARGE\bf
A De Vries-type Duality Theorem for Locally Compact Spaces -- II
}\thanks{This paper was supported by the project no. 136/2008
$``$General and Computer Topology" of the Sofia University $``$St.
Kl. Ohridski".}\\ \vspace{0.35cm}
{\large\bf Georgi D. Dimov}}

\author{}

\date{}

\begin{document}
\maketitle
\begin{abstract}
{\footnotesize
\noindent
 This paper is a continuation of \cite{D5} and in it some
 applications of the methods and results of \cite{D5} and  of
 \cite{ST, dV, R, D,D1,D2} are given. In particular,
 some generalizations of the Stone Duality Theorem \cite{ST} are obtained;   a completion theorem for
 local contact Boolean algebras is proved; a direct proof of the Ponomarev's solution \cite{P1}
 of
 Birkhoff's Problem 72 \cite{Bi} is found, and the  spaces which
 are
 co-absolute with the (zero-dimensional) Eberlein compacts are
 described.}
\end{abstract}

{\footnotesize {\em 2000 MSC:} primary 18A40, 54D45; secondary
06E15, 54C10, 54E05, 06E10.

}

\footnotetext[1]{{\footnotesize {\em E-mail address:}
gdimov@fmi.uni-sofia.bg}}

\baselineskip = \normalbaselineskip

\section*{Introduction}

This paper is a second part of the paper \cite{D5}. In it we will
use the notions, notations and results of \cite{D5} and we will
apply the methods and results obtained in \cite{D5} and in
\cite{ST,dV,R,D,D1,D2}.

In Section 1, some generalizations of the Stone Duality Theorem
\cite{ST} are obtained.
  Namely, five
  categories $\LBA$, $\ZLBA$, $\PZLBA$, $\PLBA$ and $\GBPL$ are constructed.
We show that there exists a contravariant adjunction between the
first of these categories and the category
   $\ZLC$ of zero-dimensional locally compact Hausdorff spaces (= {\em Boolean spaces})
   and continuous maps. This contravariant adjunction
   restricts to a duality between the categories $\ZLBA$ and
   $\ZLC$. The last three categories are dual to the category
   $\PZLC$ of Boolean spaces and perfect maps. The objects of the
   category $\GBPL$ are the generalized Boolean pseudolattices (= GBPLs);
the objects of   the other four categories are not GBPLs.
 Also, two subcategories $\DZLC$ and $\DPZLC$ of the category $\DLC$ dual, respectively, to the
 categories $\ZLC$ and $\PZLC$ are described.

In Section 2, we will give an explicit description of the products
in the category $\DLC$ (see \cite[Definition 2.10]{D5} for the
category $\DLC$); note that the products in the category $\DLC$
surely exist because its dual category $\HLC$ of all locally
compact Hausdorff spaces and continuous maps (see \cite[Theorem
2.14]{D5}) has sums.

In Sections 3-6, we will characterize different topological
properties of locally compact spaces by means of algebraic
characterizations of the corresponding properties of their dual
objects. As it was shown in \cite{R} by P. Roeper, the locally
compact spaces can be described (up to homeomorphism) by means of
LCAs (see \cite[Definition 1.11]{D5} for this notion), i.e. by
triples $(A,\rho,\BBBB)$. It turns out that the dual of a
topological property can have an algebraic characterization in
which only the Boolean algebra $A$ is involved. It is easy to see
that such properties are, e.g., $``$to have a given $\pi$-weight",
 $``$to have isolated points" or $``$to have a given Souslin number".
 In this paper we will study the property  $``$to have a given $\pi$-weight" and will obtain
some slight generalizations of two results of V. I. Ponomarev
\cite{P,P1}. Further on, we will characterize the dual property of
the property $``$to have a given weight" (in the class of locally
compact spaces); it is a property in whose description all three
components $A$, $\BBBB$ and $\rho$ are involved. With the help of
this characterization, we will describe the dual objects of the
metrizable locally compact spaces and we will give some easily
proved solutions of some problems analogous to Birkhoff's Problem
72 (\cite{Bi}) which was solved brilliantly by V. I. Ponomarev
\cite{P1}. We will also give a new direct solution of this
problem. Further, we will
 characterize the  spaces which are co-absolute with (zero-dimensional)
 Eberlein compacts. Let us  mention as well that there exist topological
  properties whose
dual forms are described by means of $A$ and $\BBBB$ only; such
is, for example, the property $``$to be a discrete space" (see
\cite{D2}).

 Finally, in Section 7, we will use the technique developed for the proof of our
main theorem \cite[Theorem 2.14]{D5} in order to obtain a
completion theorem for LCAs, where both the existence and the
uniqueness of the LCA-completion are proved.

For convenience of the reader, we will now repeat some of  the
notations introduced in the first part of this paper.

If $\CC$ denotes a category, we write $X\in \card\CC$ if $X$ is
 an object of $\CC$, and $f\in \CC(X,Y)$ if $f$ is a morphism of
 $\CC$ with domain $X$ and codomain $Y$.

All lattices are with top (= unit) and bottom (= zero) elements,
denoted respectively by 1 and 0. We do not require the elements
$0$ and $1$ to be distinct. We set $\2=\{0,1\}$, where $0\neq 1$.
If $(A,\le)$ is a poset and $a\in A$, we set
$\downarrow_A(a)=\{b\in A\st b\le a\}$ (we will even write
$``\downarrow(a)$" instead of $``\downarrow_A(a)$" when there is
no ambiguity); if $B\sbe A$ then we set
$\downarrow(B)=\bigcup\{\downarrow(b)\st b\in B\}$.

If $X$ is a set then we denote the power set of $X$ by $P(X)$. If
$Y$ is also a set and $f:X\lra Y$ is a function, then we will set,
for every $U\sbe X$,
 $f^\sharp(U)=\{y\in Y\st f\inv(y)\sbe U\}$.
 If
$(X,\tau)$ is a topological space and $M$ is a subset of $X$, we
denote by $\cl_{(X,\tau)}(M)$ (or simply by $\cl(M)$ or
$\cl_X(M)$) the closure of $M$ in $(X,\tau)$ and by
$\int_{(X,\tau)}(M)$ (or briefly by $\int(M)$ or $\int_X(M)$) the
interior of $M$ in $(X,\tau)$. The (positive) natural numbers are
denoted by $\mathbb{N}$ (resp., by $\mathbb{N}^+$) and the real
line -- by $\mathbb{R}$.

The  closed maps  between topological spaces are assumed to be
continuous but are not assumed to be onto. Recall that a map is
{\em perfect}\/ if it is  compact (i.e. point inverses are compact
sets) and closed. A continuous map $f:X\lra Y$ is {\em
irreducible}\/ if $f(X)=Y$ and for each proper closed subset $A$
of $X$, $f(A)\neq Y$.

For all notions and notations not defined here see \cite{D5, AHS,
J, E,  Si}.

\section{Some Generalizations of the Stone Duality Theorem}

In this section, using  Roeper's theorem \cite[Theorem 2.1]{D5},
some generalizations of the Stone Duality Theorem \cite{ST} are
obtained. A category $\LBA$ is constructed and a contravariant
adjunction between it and the category $\ZLC$ of {\em Boolean
spaces} (= zero-dimensional locally compact Hausdorff spaces) and
continuous maps is obtained. The fixed objects of this adjunction
give us a duality between the category $\ZLC$ and the subcategory
$\ZLBA$ of the category $\LBA$. Three categories $\PZLBA$, $\PLBA$
and $\GBPL$ dual to the category $\PZLC$ of Boolean spaces and
perfect maps are described. The restrictions of the obtained
duality functors  to the category $\ZHC$ of zero-dimensional
compact Hausdorff spaces (= {\em Stone spaces}) and continuous
maps  coincide with the Stone duality functor $S^t:\ZHC\lra\Bool$,
where $\Bool$ is the category of Boolean algebras and Boolean
homomorphisms. We describe as well two subcategories $\DZLC$ and
$\DPZLC$ of the category $\DLC$ which are dual, respectively,  to
the categories $\ZLC$ and $\PZLC$. Recall that complete LCAs are
abbreviated as CLCAs (see \cite[Definition 1.11]{D5}).

\begin{defi}\label{dzlc}
\rm Let $\DZLC$ (resp., $\DPZLC$) be the full subcategory of the
category $\DLC$ (resp., $\PAL$) having as objects all CLCAs
$(A,\rho,\BBBB)$ such that if $a,b\in\BBBB$ and $a\llx b$ then
there exists $c\in\BBBB$ with $c\llx c$ and $a\le c\le b$ (see
\cite[Definition 1.1]{D5} for $\llx$).
\end{defi}

\begin{theorem}\label{dzlcth}
The categories $\ZLC$ and $\DZLC$, as well as the categories
$\PZLC$ and $\DPZLC$ are dually equivalent.
\end{theorem}

\doc We will show that the contravariant functors
$\LAM^t_z=(\LAM^t)_{|\ZLC}$ and  $\LAM^a_z=(\LAM^a)_{|\DZLC}$ are
the required duality functors (see \cite[Theorem 2.14]{D5} for
$\LAM^t$ and $\LAM^a$) for the first pair of categories. Indeed,
if $X\in\card{\ZLC}$  then $\LAM^t(X)=(RC(X),\rho_X,CR(X))$ (see
\cite[ 1.3 and 1.8]{D5} for these notations) and, obviously,
$(RC(X),\rho_X,CR(X))\in\card{\DZLC}$. Conversely, if
$(A,\rho,\BBBB)\in\card{\DZLC}$ then $X=\LAM^a(A,\rho,\BBBB)$ is a
locally compact Hausdorff space. For proving that $X$ is a
zero-dimensional space, let $x\in X$ and $U$ be an open
neighborhood of $x$. Then there exist open sets $V,W$ in $X$ such
that $x\in V\sbe\cl(V)\sbe W\sbe \cl(W)\sbe U$ and $\cl(V)$,
$\cl(W)$ are compacts. Then there exist $a,b\in\BBBB$ such that
$\lag(a)=\cl(V)$ and $\lag(b)=\cl(W)$ (see \cite[(21)]{D5} for the
notation $\lag$) . Obviously, $a\llx b$. Thus, there exists
$c\in\BBBB$ such that $c\llx c$ and $a\le c\le b$. Then
$F=\lag(c)$ is a clopen subset of $X$ and $x\in F\sbe U$. So, $X$
is zero-dimensional. Now, all follows from \cite[Theorem
2.14]{D5}.

The restrictions of the obtained above duality functors to the
categories of the second pair give, according to \cite[Theorem
2.9]{D5}, the desired second duality. \sqs

\begin{defi}\label{deflba}
\rm A pair $(A,I)$, where $A$ is a Boolean algebra and $I$ is an
ideal of $A$ (possibly non proper) which is dense in $A$ (shortly,
dense ideal), is called a {\em local Boolean algebra} (abbreviated
as LBA). An LBA $(A,I)$ is called a {\em prime local Boolean
algebra} (abbreviated as PLBA) if $I=A$ or $I$ is a prime ideal of
$A$. Two LBAs $(A,I)$ and $(B,J)$ are said to be {\em isomorphic}
if there exists a Boolean isomorphism $\p:A\lra B$ such that
$\p(I)=J$.

Let $\LBA$ be the category whose objects are all LBAs and whose
morphisms are all functions $\p:(B, I)\lra(B_1, I_1)$ between the
objects of $\LBA$ such that $\p:B\lra B_1$ is a Boolean
homomorphism satisfying the following condition:

\smallskip

\noindent(LBA) For every $b\in I_1$ there exists $a\in I$ such
that $b\le \p(a)$;

\smallskip

\noindent let the composition between the morphisms of $\LBA$ be
the usual composition between functions, and the $\LBA$-identities
be the identity functions.
\end{defi}

\begin{rem}\label{remlba}
\rm Note that a prime (= maximal) ideal $I$ of a Boolean algebra
$A$ is a dense subset of $A$ iff $I$ is a non-principal ideal of
$A$. For proving this, observe first that if $I$ is a prime ideal,
$a\in A\stm\{1\}$ and $I\le a$ then $a\in I$. (Indeed, if $a\nin
I$ then $a^*\in I$ and hence $a^*\le a$, i.e. $a=1$.) Let now $I$
be dense in $A$. Suppose that $I=\downarrow(a)$ for some  $a\in
A\stm\{1\}$. Then $a^*\neq 0$. There exists $b\in I\stm\{0\}$ such
that $b\le a^*$. Since $b\le a$, we get that $b=0$, a
contradiction. Hence, $I$ is a non-principal ideal. Conversely,
let $I$ be a non-principal ideal and $b\in A\stm\{0\}$. Suppose
that $b\we a=0$, for every $a\in I$. Then $I\le b^*$. Hence
$I=\downarrow(b^*)$, a contradiction. Thus, there exists $a\in I$
such that $a\we b\neq 0$. Then $a\we b\in I\stm\{0\}$ and $a\we
b\le b$. Therefore, $I$ is a dense subset of $A$.
\end{rem}

The next obvious lemma is our motivation for  introducing the
notion of a local Boolean algebra (LBA):

\begin{lm}\label{rhoslemma}
If $(A,\rho_s,\BBBB)$ is an LCA  then $(A,\BBBB)$ is an LBA.
Conversely, for any LBA $(A,I)$, the triple $(A,\rho_s,I)$ is an
LCA. (See \cite[Example 1.2]{D5} for the notation $\rho_s$.)
\end{lm}

Since we follow Johnstone's terminology from \cite{J}, we will use
the term {\em pseudolattice} for a poset having all finite
non-empty meets and joins; the pseudolattices with a bottom will
be called {\em $\{0\}$-pseudolattices}. Recall that a distributive
$\{0\}$-pseudolattice $A$ is called a {\em generalized Boolean
pseudolattice} if it satisfies the following condition:

\smallskip

\noindent(GBPL) for every $a\in A$ and every $b,c\in A$ such that
$b\le a\le c$ there exists $x\in A$ with $a\we x=b$ and $a\vee
x=c$ (i.e., $x$ is the {\em relative complement of $a$ in the
interval} $[b,c]$).

\smallskip

Let $A$ be a distributive $\{0\}$-pseudolattice and $Idl(A)$ be
the frame of all ideals of $A$. If $J\in Idl(A)$ then we will
write $\neg_A J$ (or simply $\neg J$) for the pseudocomplement of
$J$ in $Idl(A)$ (i.e. $\neg J=\bv\{I\in Idl(A)\st I\we
J=\{0\}\}$). Note that $\neg J=\{a\in A\st (\fa b\in J)(a\we
b=0)\}$ (see Stone \cite{ST1}). Recall that an ideal $J$ of $A$ is
called {\em simple} (Stone \cite{ST1}) if $J\vee\neg J= A$. As it
is proved in \cite{ST1}, the set $Si(A)$ of all simple ideals of
$A$ is a Boolean algebra with respect to the lattice operations in
$Idl(A)$.

\begin{fact}\label{gbpapi}
(a) A distributive $\{0\}$-pseudolattice $A$ is a generalized
Boolean pseudolattice iff every principal ideal of $A$ is simple.

\smallskip

\noindent(b) If $A$ is a generalized Boolean pseudolattice then
the correspondence $e_A:A\lra Si(A)$, $a\mapsto\dar(a)$, is a
dense $\{0\}$-pseudolattice embedding of $A$ in the Boolean
algebra $Si(A)$ and the pair $(Si(A),e_A(A))$ is an LBA.

\smallskip

\noindent(c){\rm (M. Stone \cite{ST1})} An ideal of a Boolean
algebra is simple iff it is principal.
\end{fact}

\doc (a) ($\Rightarrow$)  Let $A$ be a generalized Boolean
pseudolattice and $a\in A$. We have to prove that
$\downarrow(a)\vee\neg(\downarrow(a))=A$. Let $b\in A$. Then
$c=a\we b\in [0,b]$. Hence there exists $d\in A$ such that $d\we
c=0$ and $d\vee c=b$. Thus $d\le b$, i.e. $d\we b=d$. Therefore,
$d\we a=d\we b\we a=d\we c=0$. We obtain that
$d\in\neg(\downarrow(a))$, $c\in\downarrow(a)$ and $c\vee d=b$.
So, $\downarrow(a)\vee\neg(\downarrow(a))=A$.

\smallskip

\noindent($\Leftarrow$) Let $a,b,c\in A$ and $a\in[b,c]$. Since
 $\downarrow(a)\vee\neg(\downarrow(a))=A$, we get that there
 exists $y\in \neg(\downarrow(a))$ such that $c=a\vee y$. Set
 $x=y\vee b$. Then $x\we a=(y\vee b)\we a=b\we a=b$ and $x\vee
 a=y\vee b\vee a=y\vee a=c$. So, $A$ is a generalized Boolean
 pseudolattice.

\smallskip

 \noindent(b) By (a), for every $a\in A$, $\downarrow(a)\in
 Si(A)$. Further, it is easy to see that $e_A$ is a
 $\{0\}$-pseudolattice embedding and $I=e_A(A)$ is dense in
 $Si(A)$. Let us show that $I$ is an ideal of $Si(A)$. Since $I$
 is closed under finite joins, it is enough to prove that $I$ is a
 lower set. Let $J\in Si(A)$, $a\in A$ and $J\sbe \downarrow(a)$.
 We need to show that $J$ is a principal ideal of $A$. Since $J\in
 Si(A)$, there exist $b\in J$ and $c\in\neg J$ such that $a=b\vee
 c$. We will prove that $J=\downarrow(b)$. Note first that if
 $b\ap\in J$ and $a= b\ap\vee c$ then $b=b\ap$. Indeed, we have
 that $b\ap=a\we b\ap=(b\vee c)\we b\ap=b\we b\ap$ and
 $b=a\we b=(b\ap\vee c)\we b=b\we b\ap$; thus $b=b\ap$. Let now
 $d\in J$. Then $d\le a$ and hence $a=a\vee d=(b\vee d)\vee c$.
 Since $b\vee d\in J$, we get that $b\vee d=b$, i.e. $d\le b$. So,
$J=\downarrow(b)$, and hence $J\in I$. Thus $(Si(A),e_A(A))$ is an
LBA.

\smallskip

\noindent(c) Let $B$ be a Boolean algebra and $J\in Si(B)$. Then
there exist $a\in J$ and $b\in\neg J$ such that $1=a\vee b$. Now
we obtain, as in the proof of (b), that $J=\downarrow(a)$. So,
every simple ideal of $B$ is principal. Thus, using (a), we
complete the proof.
 \sqs

\begin{nota}\label{day}
\rm Let $I$ be a proper ideal of a Boolean algebra $A$. We set
$$B_A(I)=I\cup\{a^*\st a\in I\}.$$
When there is no ambiguity, we will often write $``B(I)$" instead
of $``B_A(I)$".
\end{nota}

It is clear that $B_A(I)$ is a Boolean subalgebra of $A$ and $I$
is a prime ideal of $B_A(I)$ (see, e.g., \cite{Dw}).

\begin{fact}\label{gbplplba}
Let $(A,I)$ be an LBA. Then:

\noindent(a) $I$ is a generalized Boolean pseudolattice;

\smallskip

\noindent(b) If $(B,J)$ is a PLBA and there exists a
poset-isomorphism $\psi:J\lra I$ then $\psi$ can be uniquely
extended to a Boolean embedding $\p:B\lra A$ (and $\p(B)=B_A(I)$);
in particular, if $(A,I)$ is also a PLBA then $\p$ is a Boolean
isomorphism and an isomorphism between LBAs $(A,I)$ and $(B,J)$;

\smallskip

\noindent(c) There exists a bijective correspondence between the
class of all (up to isomorphism) generalized Boolean
pseudolattices and the class of all (up to isomorphism) PLBAs.
\end{fact}

\doc (a) Obviously, for every $a\in I$,
$\neg_I(\dar(a))=I\cap\dar_A(a^*)$; then, clearly,
$\dar(a)\vee\neg_I(\dar(a))=I$. Now apply \ref{gbpapi}(a).

\smallskip

\noindent(b) By \cite[Theorem 12.5]{Si}, $\psi$ can be uniquely
extended to a Boolean isomorphism $\psi\ap:B\lra B_A(I)$. Now,
define $\p:B\lra A$ by $\p(b)=\psi\ap(b)$, for every $b\in B$.

\smallskip

\noindent(c) For every PLBA $(A,I)$, set $f(A,I)=I$. Then, by (a),
$I$ is a generalized Boolean pseudolattice. Conversely, if $I$ is
a generalized Boolean pseudolattice then there exists a dense
embedding $e:I\lra Si(I)$ (see Fact \ref{gbpapi}(b)). Thus,
setting $g(I)=(B_{Si(I)}(e(I)),e(I))$, we get that $g(I)$ is a
PLBA. Now, using (b), we obtain that for every PLBA $(A,I)$,
$g(f(A,I))$ is isomorphic to $(A,I)$. Finally, it is clear that
for every generalized Boolean pseudolattice $I$, $f(g(I))$ is
isomorphic to $I$. \sqs

\begin{lm}\label{crhoslemma}
Let $(A,I)$ be an LBA and $\s\sbe A$. Then $\s$ is a bounded
cluster in $(A,\rho_s,I)$ iff it
 is a bounded ultrafilter in $A$ (see \cite[Definition 1.15]{D5} for the last notion).
\end{lm}

\doc Let $C=C_{\rho_s}$ be the Alexandroff extension of the
relation $\rho_s$ relatively to the LCA $(A,\rho_s,I)$ (see
\cite[Definition 1.13]{D5} for $C_{\rho_s}$ and \ref{rhoslemma}
for $(A,\rho_s,I)$).

Using \cite[Theorem 1.8]{D5} and \cite[Corollary 1.9]{D5}, we
obtain that:
 [$\s\sbe A$ is a bounded cluster in $(A,\rho_s,I)$] $\iff$ [$\s$ is a
cluster in $(A,C)$ and $\s\cap I\nes$] $\iff$ [there exists a
bounded ultrafilter $u$ in $A$ such that $\s=\s_u$]. Hence
$\s=\{a\in A\st (\fa b\in u)(a C_{\rho_s} b)\}$. Note that $u\cap
I$ is a filter base of $u$. (Indeed, since $u$ is bounded, there
exists $a_0\in u\cap I$; then, for every $a\in u$, $b=a\we a_0\in
u\cap I$ and $b\le a$.) Thus $\s=\{a\in A\st (\fa b\in u\cap I)(a
C_{\rho_s} b)\}=\{a\in A\st (\fa b\in u\cap I)(a\we b\neq
0)\}=\{a\in A\st (\fa b\in u)(a\we b\neq 0)\}=u$. \sqs

\begin{notas}\label{kxckx}
\rm Let $X$ be a topological space. We will denote by $CO(X)$ the
set of all clopen subsets of $X$,
and by $CK(X)$ the set of all clopen compact subsets of $X$. For
every $x\in X$, we set
$$  u_x^{CO(X)}=\{F\in CO(X)\st x\in F\}
.$$
When there is no ambiguity, we will write $``u_x^C$" instead of
$``u_x^{CO(X)}$".
\end{notas}

Recall that a {\em contravariant adjunction}\/ between two
categories $\AA$ and $\BB$ consists of two contravariant functors
$T:\AA\lra\BB$ and $S:\BB\lra \AA$ and two natural transformations
$\eta:Id_{\BB}\lra T\circ S$ and $\ep:Id_{\AA}\lra S\circ T$ such
that $T(\ep_A)\circ\eta_{TA}=id_{TA}$ and $S(\eta_B)\circ
\ep_{SB}=id_{SB}$, for all $A\in\card{\AA}$ and $B\in\card{\BB}$
(here, as usual,  $Id$ is the identity functor and $id$ is the
identity morphism). The pair $(S,T)$ is a duality iff $\eta$ and
$\ep$ are natural isomorphisms.

\begin{theorem}\label{genstonec}
There exists a contravariant adjunction between the category\/
$\LBA$ and the category\/ $\ZLC$ of  locally compact
zero-dimensio\-nal Hausdorff spaces and
 continuous maps.
\end{theorem}

\doc  We will first define two
contravariant functors $\TE^a:\LBA\lra\ZLC$ and
$\TE^t:\ZLC\lra\LBA$.

Let $X\in\card{\ZLC}$. Define $$\TE^t(X)=(CO(X),  CK(X)).$$
Obviously, $\TE^t(X)$ is an LBA.

Let $f\in\ZLC(X,Y)$. Define $\TE^t(f):\TE^t(Y)\lra\TE^t(X)$ by the
formula
\begin{equation}\label{deftetaf}
\TE^t(f)(G)=f\inv(G), \ \  \fa G\in CO(Y).
\end{equation}
Set $\p_f=\TE^t(f)$. Clearly, $\p_f$ is a Boolean homomorphism
between $CO(Y)$ and $CO(X)$. If $F\in CK(X)$ then $f(F)$ is a
compact subset of $Y$. Since $CK(Y)$ is an open base of the space
$Y$ and $CK(Y)$ is closed under finite unions, we get that there
exists $G\in CK(Y)$ such that $f(F)\sbe G$. Then $F\sbe
f\inv(G)=\p_f(G)$. So, $\p_f$ satisfies condition (LBA). Therefore
$\p_f$ is a $\LBA$-morphism, i.e. $\TE^t(f)$ is well-defined.

Now we get easily that $\TE^t$ is a contravariant functor.

For every LBA $(B, I)$, set $$\TE^a(B, I)=\Psi^a(B,\rho_s,I)$$
(see \cite[(13) and (15)]{D5} for $\Psi^a$ and \ref{rhoslemma} for
the fact that $(B,\rho_s,I)$ is an LCA). Then \cite[Theorem
2.1(a)]{D5} implies that $L=\TE^a(B, I)$ is a locally compact
Hausdorff space. Since for any $a\in B$ we have that
$a\ll_{\rho_s} a$, we get that $\lbg  (B)\sbe CO(L)$. By
\cite[(24)]{D5}, $\lbg  (I)$ is an open base of $L$. Thus, $L$ is
a zero-dimensional space. So, $\TE^a(B,I)\in\card{\ZLC}$.

Let $\p\in\LBA((B, I),(B_1, J))$. We define the map
$$\TE^a(\p):\TE^a(B_1, J)\lra\TE^a(B, I)$$ by the formula
\begin{equation}\label{deftetaphi}
\TE^a(\p)(u\ap)=\p\inv(u\ap), \ \  \fa u\ap\in\TE^a(B_1,J).
\end{equation}
Set $f_\p=\TE^a(\p)$, $L=\TE^a(B, I)$ and $M=\TE^a(B_1, J)$.

By Lemma \ref{crhoslemma}, \cite[(14)]{D5} and \cite[(15)]{D5}, if
$(B\ap,I\ap)$ is a LBA then the set $\TE^a(B\ap,I\ap)$ consists of
all bounded ultrafilters of $B\ap$ (i.e., those ultrafilters $u$
of $B\ap$ for which $u\cap I\ap\nes$). Since any $\LBA$-morphism
is a Boolean homomorphism, we get that the inverse image of an
ultrafilter is an ultrafilter.

So, let $u\ap\in M$. Then $u\ap$ is a bounded ultrafilter in
$B_1$. Set $u=f_\p(u\ap)$. Then, as we have seen, $u$ is an
ultrafilter in $B$. We have to show that $u$ is bounded. Indeed,
since $u\ap$ is bounded, there exists $b\in u\ap\cap J$. By (LBA),
there exists $a\in I$ such that $\p(a)\ge b$. Then $\p(a)\in
u\ap$, and hence, $a\in u$. Thus $a\in u\cap I$. Therefore,
$f_\p:M\lra L$.

We will show that $f_\p$ is a continuous function. Let $u\ap\in M$
and $u=f_\p(u\ap)$. Let $a\in B$ and $u\in\lbg (a)(=\int(\lbg
(a)))$. Then $a\in u$. Hence $\p(a)\in u\ap$, i.e.
$u\ap\in\l_{B_1}^g(\p(a))$. We will prove that
\begin{equation}\label{contfphi}
f_\p(\l_{B_1}^g(\p(a)))\sbe\lbg  (a).
\end{equation}
Indeed, let $v\ap\in\l_{B_1}^g(\p(a))$. Then $\p(a)\in v\ap$. Thus
$a\in f_\p(v\ap)$, i.e. $f_\p(v\ap)\in\lbg  (a)$. So,
(\ref{contfphi}) is proved. Since $\{\lbg  (a)\st a\in B\}$ is an
open base of $L$, we get that $f_\p$ is a continuous function. So,
$$\TE^a(\p)\in\ZLC(\TE^a(B_1, J),\TE^a(B, I)).$$ Now it becomes
obvious that $\TE^a$ is a contravariant functor.

Let $X\in\card{\ZLC}$. Then it is easy to see that for every $x\in
X$, $ u_x^C$ is an ultrafilter in $CO(X)$ and hence, by Lemma
\ref{crhoslemma} and the fact that $ u_x^C$ contains always
elements of $CK(X)$,  we get that $$ u_x^C\in\TE^a(CO(X),
CK(X)).$$ We will show that the map $t_X^C:X\lra\TE^a(\TE^t(X))$
defined by $t_X^C(x)= u_x^C$, for every $x\in X$, is a
homeomorphism. Set $L=\TE^a(\TE^t(X))$ and $B=CO(X)$, $I=CK(X)$.
We will prove that $t_X^C$ is a continuous map. Let $x\in X$,
$F\in I$  and $ u_x^C\in\lbg  (F)$. Then $F\in\scx$ and hence,
$x\in F$. It is enough to show that $t_X^C(F)\sbe\lbg (F)$. Let
$y\in F$. Then $F\in\scy=\tcx(y)$. Hence $\tcx(y)\in\lbg (F)$. So,
$t_X^C(F)\sbe\lbg  (F)$. Since $\lbg  (I)$ is an open base of $L$,
we get that $\tcx$ is a continuous map. Let us show that $\tcx$ is
a bijection.  Let $u\in L$. Then $u$ is a bounded ultrafilter in
$(B,\rho_s,I)$. Hence, there exists $F\in u\cap I$. Since $F$ is
compact, we get that $\bigcap u\nes$. Suppose that $x,y\in\bigcap
u$ and $x\neq y$. Then there exist $F_x,F_y\in I$ such that $x\in
F_x$, $y\in F_y$ and $F_x\cap F_y=\ems$. Since, clearly,
$F_x,F_y\in u$, we get a contradiction. So, $\bigcap u=\{x\}$ for
some $x\in X$. It is clear now that $u=\scx$, i.e., $u=\tcx(x)$
and  $u\neq\tcx(y)$, for $y\in X\stm\{x\}$. So, $\tcx$ is a
bijection. For showing that $(\tcx)\inv$ is a continuous function,
let $\scx\in L$. Then $(\tcx)\inv(\scx)=x$. Let $F\in I$ and $x\in
F$. Then $F\in\scx$ and thus $\scx\in\lbg (F)$. We will prove that
$(\tcx)\inv(\lbg  (F))\sbe F$. Since $I$ is a base of $X$, this
will imply that $(\tcx)\inv$ is a continuous function. So, let
$y\in(\tcx)\inv(\lbg  (F))$. Then $\tcx(y)\in\lbg  (F)$, i.e.
$F\in\scy$. Then $y\in F$. Therefore, $\tcx$ is a homeomorphism.

We will show that
%
$$t^{C}:Id_{\ZLC}\lra   \TE^a\circ  \TE^t,$$
%
 defined by
%
 $t^{C}(X)=\tcx, \ \ \fa X\in\card\ZLC,$
%
is a natural isomorphism.

Let $f\in\ZLC(X,Y)$ and $\fs= \TE^a(\TE^t(f))$. We have to show
that $\fs\circ \tcx=\tcy\circ f$. Let $x\in X$. Then
$\fs(\tcx(x))=\fs( u_x^{CO(X)})$ and $(\tcy\circ
f)(x)=u_{f(x)}^{CO(Y)}$. Set $y=f(x)$, $u_x= u_x^{CO(X)}$ and
$u_y=u_{f(x)}^{CO(Y)}$. We will prove that
$$\fs(u_x)=u_y.$$
Let $\p=\TE^t(f)$. Then $\fs=\TE^a(\p)(=f_\p)$. Hence,
$\fs(u_x)=\p\inv(u_x)=\{G\in CO(Y)\st \p(G)\in u_x\}= \{G\in
CO(Y)\st x\in\p(G)\}=\{G\in CO(Y)\st x\in f\inv(G)\}=\{G\in
CO(Y)\st f(x)\in G\}=u_y$. So, $t^C$ is a natural isomorphism.

Let $(B, I)$ be an LBA and $L=\TE^a(B, I)$. Then, by
\cite[(22)]{D5}, $\lbg:(B,\rho_s,I)\lra(RC(L),\rho_L,CR(L))$ is a
dense LCA-embedding. Also, obviously, $\lbg(B)\sbe CO(L)$ and
$\lbg(I)\sbe CK(L)$. We denote by $\l_{(B,I)}^C$ the restriction
$\l_{(B,I)}^C:(B,I)\lra(CO(L),CK(L))$ of $\lbg$, i.e.
$\l_{(B,I)}^C(b)=\lbg(b)$, for every $b\in B$;  we will write
sometimes $``\lbc$" instead of $``\l_{(B,I)}^C$". Note that
$\l_{(B,I)}^C:(B,I)\lra\TE^t(\TE^a(B,I))$. We will prove that
$$\l^C: Id_{\LBA}\lra  \TE^t\circ  \TE^a,\mbox{ where } \l^C(B,
I)=\lbc, \ \ \fa (B, I)\in\card\LBA,$$
%
 is a natural
transformation.

Let $\p\in\LBA((B, I),(B_1, J))$ and $\ps=\TE^t(\TE^a(\p))$. We
have to prove that $\l_{B_1}^C\circ\p=\ps\circ \l_B^C$. Set $f=
\TE^a(\p)$ and $M=\TE^a(B_1, J)$. Then $\ps=\TE^t(f)(=\p_f)$. Let
$a\in B$. Then $\ps(\lbc(a))=f\inv(\lbc(a))=\{u\in M\st
f(u)\in\lbc(a)\}=\{u\in M\st a\in f(u)\}=\{u\in M\st
a\in\p\inv(u)\}=\{u\in M\st \p(a)\in u\}=\l_{B_1}^C(\p(a))$. So,
$\l^C$ is a natural transformation.

Let us show that $\TE^t(t_X^C)\circ\l^C_{\TE^t(X)}=id_{\TE^t(X)}$,
for every $X\in\card{\ZLC}$. Indeed, let $X\in\card{\ZLC}$ and
$Y=\TE^a(\TE^t(X))$. Then $\TE^t(t_X^C):\TE^t(Y)\lra\TE^t(X)$,
$G\mapsto (t_X^C)\inv(G)$ for every $G\in\TE^t(Y)=(CO(Y),CK(Y))$.
Let $F\in CO(X)$. Then
$(\TE^t(t_X^C)\circ\l^C_{\TE^t(X)})(F)=(t_X^C)\inv(\l^C_{\TE^t(X)}(F))=H$.
We have to show that $F=H$. Since $t_X^C(H)=\l^C_{\TE^t(X)}(F)$,
we get that $\{u_x^C\st x\in H\}=\{u\in Y\st F\in u\}$. Thus $x\in
H\iff F\in u_x^C\iff x\in F$. Therefore, $F=H$.

Finally, we will prove that $\TE^a(\l_{(A,I)}^C)\circ
t^C_{\TE^a(A,I)}=id_{\TE^a(A,I)}$ for every $(A,I)\in\card{\LBA}$.
So, let $(A,I)\in\card{\LBA}$ and $X=\TE^a(A,I)$. We have that
$f=\TE^a(\l_{(A,I)}^C):\TE^a(CO(X),CK(X))\lra X$ is defined by
$u\mapsto(\l_{(A,I)}^C)\inv(u)$, for every bounded ultrafilter $u$
in $(CO(X),CK(X))$. Let $x\in X$. Then
$f(t_X^C(x))=f(u_x^C)=(\l_{(A,I)}^C)\inv(u_x^C)=y$.
 We have to show that $x=y$. Indeed, for
every $a\in A$, we get that $a\in y\iff
a\in(\l_{(A,I)}^C)\inv(u_x^C)\iff  \l_{(A,I)}^C(a)\in u_x^C\iff
x\in\l_{(A,I)}^C(a)\iff a\in x$. Therefore, $x=y$.

We have proved that $(\TE^t,\TE^a,\l^C,t^C)$ is a contravariant
adjunction between the categories $\ZLC$ and $\LBA$. Moreover, we
have shown that $t^C$ is even a natural isomorphism.
 \sqs

\begin{defi}\label{defzlba}
\rm An LBA $(B, I)$  is called a {\em ZLB-algebra} (briefly, {\em
ZLBA}) if, for every $J\in Si(I)$, the join $\bv_B J$($=\bv_B
\{a\st a\in J\}$) exists.

Let $\ZLBA$ be the full subcategory of the category $\LBA$ having
as objects  all ZLBAs.
\end{defi}

\begin{exa}\label{zlbaexa}
\rm Let $B$ be a Boolean algebra. Then the pair $(B,B)$ is a ZLBA.
This follows from Fact \ref{gbpapi}(c).
\end{exa}

\begin{rem}\label{remzlba}
\rm  Note  that if $A$ and $B$ are Boolean algebras then any
Boolean homomorphism $\p:A\lra B$ is a $\ZLBA$-morphism between
the ZLBAs $(A,A)$ and $(B,B)$. Hence, the full subcategory $\B$ of
the category $\ZLBA$ whose objects are all ZLBAs of the form
$(A,A)$ is isomorphic (it can be even said that it coincides) with
the category $\Bool$ of Boolean algebras and Boolean
homomorphisms.
\end{rem}

We will need the following result of M. Stone \cite{ST}:

\begin{pro}\label{stcox}{\rm (M. Stone
\cite[Theorem 5(3)]{ST})} Let $X\in\card{\ZLC}$. Then the map
$\SI:Si(CK(X))\lra CO(X)$, $J\mapsto \bv_{RC(X)}J$, is a Boolean
isomorphism.
\end{pro}

\doc For completeness of our exposition, we will verify this fact.
Let $J\in Si(CK(X))$. Set $U=\bigcup\{F\st F\in J\}$ and
$V=\bigcup\{G\st G\in\neg J\}$. Obviously, $U$ and $V$ are
disjoint open subsets of $X$. We will show that $U\cup V=X$.
Indeed, let $x\in X$. Then there exists $H\in CK(X)$ such that
$x\in H$. Since $J\vee\neg J=CK(X)$, we get that there exist $F\in
J$ and $G\in\neg J$ such that $H=F\cup G$. Thus $x\in F$ or $x\in
G$, and hence, $x\in U$ or $x\in V$. So, $U$ is a clopen subset of
$X$. Thus $U\in CO(X)$ and $U=\bv_{RC(X)}J=\bv_{CO(X)}J$.
Conversely, it is easy to see that if $U\in CO(X)$ then $J=\{F\in
CK(X)\st F\sbe U\}\in Si(CK(X))$. This implies easily that $\SI$
is a Boolean isomorphism. \sqs

\begin{pro}\label{cox}
Let $(B,I)$ be an LBA and $X=\TE^a(B,I)$. Then:

\smallskip

\noindent(a) $\lbg(I)=CK(L)$;

\smallskip

 \noindent(b) $(B,I)$ is a ZLBA iff $\lbg(B)=CO(X)$.
\end{pro}

\doc (a) We have that $\lbg(B)\sbe CO(X)$ and hence,
$\lbg(I)\sbe CO(X)\cap CR(X)=CK(X)$. Conversely, if $F\in CK(X)$
then the facts that $\lbg(I)$ is an open base of $X$ and $\lbg(I)$
is closed under finite unions imply that $F\in\lbg(I)$. Thus,
$\lbg(I)=CK(X)$.

\smallskip

 \noindent(b) Let $(B,I)$ be a ZLBA. We will prove that $\lbg(B)=CO(X)$. Let
$U\in CO(X)$ and $J\ap=\{F\in CK(X)\st F\sbe U\}$. Then $J\ap$ is
a simple ideal of $CK(X)$ and $\bv_{RC(X)}J\ap=U$. Since the
restriction $\p:I \lra CK(X)$ of $\lbg$ is a $\{0\}$-pseudolattice
isomorphism, we get that $J=\p\inv(J\ap)$ is a simple ideal of
$I$. Set $b_J=\bv_B J$ and $C=\lbg(B)$ (note that the join $\bv_B
J$ exists  because $(B,I)$ is  a ZLBA). Now, the restriction
$\psi:B\lra C$ of $\lbg$ is a Boolean isomorphism, and hence
$\lbg(b_J)=\psi(b_J)=\psi(\bv_B J)=\bv_C\psi(J)=\bv_C J\ap$. The
fact that $C$ is a dense Boolean subalgebra of the Boolean algebra
$RC(X)$ implies that $C$ is a regular subalgebra of $RC(X)$. Thus
$\bv_C J\ap=\bv_{RC(X)} J\ap=U$. Therefore, $\lbg(b_J)=U$. So, we
have proved that $\lbg(B)=CO(X)$.

Let now $(B,I)$ be an LBA and $\lbg(B)=CO(X)$. Then, as above, the
restriction $\psi:B\lra CO(X)$ of $\lbg$ is a Boolean isomorphism.
Let $J\in Si(I)$. Since, by (a), the restriction of $\psi$ to $I$
is a 0-pseudolattice isomorphism between $I$ and $CK(X)$, we get
that $\psi(J)\in Si(CK(X))$. Then, by \ref{stcox},
$U=\bigcup\{F\st F\in \psi(J)\}(=\bigcup\{\psi(a)\st a\in J\})$ is
a clopen subset of $X$. Therefore, the join
$\bigvee_{CO(X)}\{\psi(a)\st a\in J\}$ exists. Since
$\psi\inv:CO(X)\lra B$ is a Boolean isomorphism, we obtain that
$\psi\inv(U)=\psi\inv(\bigvee_{CO(X)}\{\psi(a)\st a\in
J\})=\bigvee_B\{\psi\inv(\psi(a))\st a\in J\}=\bigvee_B\{a\st a\in
J\}$. Hence, the join $\bigvee_B J$ exists. Thus, $(B,I)$ is a
ZLBA.
 \sqs


\begin{theorem}\label{genstonecnew}
The categories\/ $\ZLC$ and\/ $\ZLBA$ are  dually equivalent.
\end{theorem}

\doc In Theorem
\ref{genstonec}, we constructed a contravariant adjunction
$$(\TE^t,\TE^a,\l^C,t^C)$$ between the categories $\ZLC$ and $\LBA$,
where $t^C$ was even a natural isomorphism. Let us check that the
functor $\TE^t$ is in fact a functor from the category $\ZLC$ to
the category $\ZLBA$. Indeed, let $X\in\card{\ZLC}$. Then
$\TE^t(X)=(CK(X),CO(X))$. As it follows from \ref{stcox}, for
every $J\in Si(CK(X))$, $\bigvee_{CO(X)} J$ exists. Hence,
$\TE^t(X)\in\card{\ZLBA}$. So, the restriction
 $$\TE^t_d:\ZLC\lra\ZLBA$$ of the contravariant functor $\TE^t:\ZLC\lra\LBA$ is well-defined. Further, by
Proposition \ref{cox}, the natural transformation $\l^C$ becomes a
natural isomorphism exactly on the subcategory $\ZLBA$ of the
category $\LBA$. We will denote by $$\TE^a_d:\ZLBA\lra\ZLC$$ the
restriction of the contravariant functor $\TE^a$ to the category
$\ZLBA$. All this shows that there is a duality between the
categories $\ZLC$ and $\ZLBA$. \sqs


\begin{cor}\label{stoneth}{\rm (Stone Duality Theorem \cite{ST})}
The categories $\Bool$ and $\ZHC$ are dually equivalent.
\end{cor}

\doc Obviously, the restriction of the contravariant functor $\TE^a_d$ to
the subcategory $\B$ of the category $\ZLBA$ (see \ref{remzlba}
for the notation $\B$) produces a duality between the categories
$\B$ and $\ZHC$. \sqs

\begin{cor}\label{stfact}
For every ZLBA $(B,I)$, the map $\SI_{(B,I)}:Si(I)\lra B$, where
$\SI_{(B,I)}(J)=\bv_B \{a\st a\in J\}$ for every $J\in Si(I)$, is
 a Boolean isomorphism.
\end{cor}

\doc Let $L=\TE^a_d(B,I)$ (see the proof of Theorem \ref{genstonecnew}
for the notation $\TE^a_d$). Then, as it was shown in the proof of
Theorem \ref{genstonecnew}, the map $\l_B^C:(B,I)\lra
(CO(L),CK(L))$, where $\l_B^C(b)=\lbg(b)$ for every $b\in B$, is a
$\ZLBA$-isomorphism. By \ref{stcox}, the map
$\SI=\SI_{(CO(L),CK(L))}:Si(CK(L))\lra CO(L)$, $J\mapsto
\bv_{CO(L)} J$, is a Boolean isomorphism. Define a map
$\l_B\ap:Si(I)\lra Si(CK(L))$ by the formula
$\l_B\ap(J)=\l_B^C(J)$, for every $J\in Si(I)$. Then, obviously,
$\l_B\ap$ is a Boolean isomorphism and
$\SI_{(B,I)}=(\l_B^C)\inv\circ\SI\circ\l_B\ap$. Thus $\SI_{(B,I)}$
is a Boolean isomorphism. \sqs

\begin{defi}\label{defplbap}
\rm Let $\PZLBA$ be the subcategory of the category $\ZLBA$,
having the same objects (i.e. $\card{\PZLBA}=\card{\ZLBA}$), whose
morphisms $\p:(A,I)\lra (B,J)$ satisfy the following additional
condition:

\smallskip

\noindent(PLBA) $\p(I)\sbe J$.
\end{defi}

\begin{theorem}\label{genstonep}
The category\/ $\PZLC$ of all locally compact zero-dimensio\-nal
Hausdorff spaces and all perfect maps between them is dually
equivalent to the category\/ $\PZLBA$.
\end{theorem}

\doc Let $f\in\PZLC(X,Y))$. Then, as we have seen in the proof of Theorem \ref{genstonec},
 $\TE^t_d(f):\TE^t_d(Y)\lra\TE^t_d(X)$ is defined by
the formula $\TE^t_d(f)(G)=f\inv(G), \fa G\in CO(Y)$. Set
$\p_f=\TE^t_d(f)$. Since $f$ is a perfect map, we have that for
any $K\in CK(Y)$, $\p_f(K)=f\inv(K)\in CK(X)$. Hence, $\p_f$
satisfies condition (PLBA). Thus, $\p_f$ is a $\PZLBA$-morphism.
So, the restriction $\TE^t_p$ of the duality functor  $\TE^t_d$ to
the subcategory $\PZLC$ of the category $\ZLC$ is a contravariant
functor from $\PZLC$ to $\PZLBA$.

Let $\p\in\PZLBA((A, I),(B, J))$. Then the map
$\TE^a_d(\p):\TE^a_d(B, J)\lra\TE^a_d(A, I)$ was defined in
Theorem \ref{genstonec} by the formula
$\TE^a_d(\p)(u\ap)=\p\inv(u\ap)$, $\fa u\ap\in\TE^a_d(B,J)$. Set
$f_\p=\TE^a_d(\p)$, $L=\TE^a_d(A, I)$ and $M=\TE^a_d(B, J)$.

Let $a\in I$. We will show that $f_\p\inv(\lag(a))$ is compact. We
have, by (PLBA), that $\p(a)\in J$. Let us prove that
\begin{equation}\label{perfphi}
\lbg(\p(a))=f_\p\inv(\lag(a)).
\end{equation}
  Let $u\ap\in f_\p\inv(\lag(a))$. Then
$u=f_\p(u\ap)\in\lag(a)$, i.e. $a\in u$. Thus $\p(a)\in u\ap$, and
hence $u\ap\in\lbg(\p(a))$. Therefore, $\lbg(\p(a))\spe
f_\p\inv(\lag(a))$. Now, (\ref{contfphi}) implies that
$\lbg(\p(a))=f_\p\inv(\lag(a))$. Since $\lbg(\p(a))$ is compact,
we get that $f_\p\inv(\lag(a))$ is compact. Let now $K$ be a
compact subset of $L$. Since $\lag(I)$ is an open base of $L$ and
$\lag(I)$ is closed under finite unions, we get that there exists
$a\in I$ such that $K\sbe\lag(a)$. Then $f_\p\inv(K)\sbe
f_\p\inv(\lag(a))$, and hence, as a closed subset of a compact
set, $f_\p\inv(K)$ is compact. This implies that $f_\p$ is a
perfect map (see, e.g.,\cite{E}). Therefore,  the restriction
$\TE^a_p$ of the duality functor  $\TE^a_d$ to the subcategory
$\PZLBA$ of the category $\ZLBA$ is a contravariant functor from
$\PZLBA$ to $\PZLC$. The rest follows from Theorem
\ref{genstonecnew}. \sqs

The above theorem can be stated in a better form. We will do this
now.

\begin{defi}\label{defplba}
\rm Let $\PLBA$ be  the subcategory of the category $\LBA$ whose
objects are all PLBAs and whose morphisms are all $\LBA$-morphisms
$\p:(A,I)\lra(B,J)$ between the objects of $\PLBA$  satisfying
condition (PLBA).
\end{defi}

\begin{rem}\label{remplba}
\rm It is obvious that $\PLBA$ is indeed a category. Note also
that any Boolean homomorphism $\p:A\lra B$ is a $\PLBA$-morphism
between the PLBAs $(A,A)$ and $(B,B)$. Hence, the full subcategory
$\B$ of the category $\PLBA$ whose objects are all PLBAs of the
form $(A,A)$ is isomorphic (it can be even said that it coincides)
with the category $\Bool$ of Boolean algebras and Boolean
homomorphisms.
\end{rem}

\begin{theorem}\label{genstone}
The category\/ $\PZLC$  is dually equivalent to the category\/
$\PLBA$.
\end{theorem}

\doc In virtue of Theorem \ref{genstonep}, it is enough to show
that the categories $\PLBA$ and $\PZLBA$ are equivalent.

Let $(B,I)$ be a ZLBA. Set $A=B_B(I)$ (see \ref{day} for the
notations). Then, obviously, $(A,I)$ is a PLBA. Set
$E^z(B,I)=(A,I)$.

If $\p\in\PZLBA((B_1,I_1),(B_2,I_2))$ then let $E^z(\p)$ be the
restriction of  $\p$ to $E^z(B_1,I_1)$. Then, clearly,
$E^z(\p)\in\PLBA(E^z(B_1,I_1),E^z(B_2,I_2))$. It is evident that
$E^z$ is a (covariant) functor from $\PZLBA$ to $\PLBA$.

Let $(A,I)$ be a PLBA. Then, by \ref{gbplplba}(a), $I$ is a
generalized Boolean pseudolattice. Hence, according to
\ref{gbpapi}(b), the map $e_I:I\lra Si(I)$, where
$e_I(a)=\downarrow(a)$, is a dense embedding of $I$ in the Boolean
algebra $Si(I)$ and the pair $(Si(I),e_I(I))$ is an LBA.   Set
$I\ap=e_I(I)$ and $E^p(A,I)=(Si(I),I\ap)$. Then, for every $J\in
Si(I)$, $\bv_{Si(I)} e_I(J) =\bv_{Si(I)}\{\downarrow(a)\st a\in
J\}=J$. This implies that $(Si(I),I\ap)\in\card{\PZLBA}$.

Let $\p\in\PLBA((A_1,I_1),(A_2,I_2))$. Let the map $\p\ap=E^p(\p)$
be defined by the formula
$\p\ap(J_1)=\bigcup\{\downarrow(\p(a))\st a\in J_1\}$, for every
$J_1\in Si(I_1)$. We will show that $\p\ap$ is a $\PZLBA$-morphism
between $E^p(A_1,I_1)$ and $E^p(A_2,I_2)$. Obviously,
$\p\ap(\{0\})=\{0\}$ and, thanks to conditions (LBA) and (PLBA),
$\p\ap(I_1)=I_2$. Let $J_1\in Si(I_1)$. Set $J_2=\p\ap(J_1)$. Then
condition (PLBA) and the fact that $\p$ is a homomorphism imply
that $J_2$ is an ideal of $I_2$. Let us show that $J_2\vee\neg
J_2=I_2$. Indeed, let $a_2\in I_2$. Then condition (LBA) implies
that there exists $a_1\in I_1$ such that $a_2\le \p(a_1)$. Since
$J_1\vee\neg J_1=I_1$, there exist $a_1\ap\in J_1$ and
$a^{\prime\prime}_1\in\neg J_1$ such that $a_1=a_1\ap\vee
a_1^{\prime\prime}$. Then $a_2=(\p(a_1\ap)\we
a_2)\vee(\p(a_1^{\prime\prime})\we a_2)$. Obviously,
$(\p(a_1\ap)\we a_2)\in J_2$. We will prove that
$(\p(a_1^{\prime\prime})\we a_2)\in\neg J_2$. It is enough to show
that $\p(a_1^{\prime\prime})\in\neg J_2$. Let $b_2\in J_2$. Then,
by the definition of $J_2$, there exists $b_1\in J_1$ such that
$b_2\le\p(b_1)$. Since $b_1\we a_1^{\prime\prime}=0$, we get that
$\p(b_1)\we \p(a_1^{\prime\prime})=0$. Thus
$\p(a_1^{\prime\prime})\we b_2=0$. Therefore,
$\p(a_1^{\prime\prime})\in\neg J_2$. So, $J_2\in Si(I_2)$. Note
that this implies that
$\p\ap(J_1)=\bv_{Si(I_2)}\{\downarrow(\p(a))\st a\in J_1\}$. The
above arguments show also that $\p\ap(\neg J_1)\sbe
\neg\p\ap(J_1)$, for every $J_1\in Si(I_1)$. In fact, there is an
equality here, i.e. $\p\ap(\neg J_1)= \neg\p\ap(J_1)$. Indeed, let
$b_2\in\neg\p\ap(J_1)$. Then $b_2\we a_2=0$, for every
$a_2\in\p\ap(J_1)$. By condition (LBA), there exists $a_1\in I_1$
such that $b_2\le\p(a_1)$. We have again that there exist
$a_1\ap\in J_1$ and $a^{\prime\prime}_1\in\neg J_1$ such that
$a_1=a_1\ap\vee a_1^{\prime\prime}$. Then $b_2=(\p(a_1\ap)\we
b_2)\vee(\p(a_1^{\prime\prime})\we b_2)=\p(a_1^{\prime\prime})\we
b_2$. Thus, $b_2\le\p(a_1^{\prime\prime})$. This shows that
$b_2\in\p\ap(\neg J_1)$. Further, if $J, J\ap\in Si(I_1)$ then
$\p\ap(J)\we\p\ap(J\ap)=\p\ap(J)\cap\p\ap(J\ap)=\bigcup
\{\downarrow(a)\we\downarrow(b)\st a\in J, b\in J\ap\}=
\bigcup\{\downarrow(a)\st a\in J\cap J\ap\}=\p\ap(J\cap
J\ap)=\p\ap(J\we J\ap)$. Therefore, $\p\ap:Si(I_1)\lra Si(I_2)$ is
a Boolean homomorphism. Since, for every $a\in I_1$,
$\p\ap(\downarrow(a))=\downarrow(\p(a))$, we have that
$e_{I_2}\circ\p_{\mid I_1}=\p\ap\circ e_{I_1}$. This shows that
$\p\ap\in\PZLBA(E^p(A_1,I_1),E^p(A_2,I_2))$. Now one can easily
see that $E^p$ is a (covariant) functor between the categories
$\PLBA$ and $\PZLBA$.

Finally, we have to verify that the compositions $E^p\circ E^z$
and $E^z\circ E^p$ are naturally isomorphic to the corresponding
identity functors.

Let us start with the composition $E^z\circ E^p$.

Let $(A,I)$ be a PLBA. Then, as we have seen above, the map
$e_I:I\lra Si(I)$, where $e_I(a)=\downarrow(a)$, is a dense
embedding of $I$ in the Boolean algebra $Si(I)$ and the pair
$(Si(I),e_I(I))$ is an LBA. Now \ref{gbplplba}(b) implies that the
map $(e_I)_{\upharpoonright I}: I\lra e_I(I)$ can be extended to a
Boolean isomorphism $e_{(A,I)}:A\lra B_{Si(I)}(e_I(I))$. (Note
that $A=I\cup I^*$ and $B_{Si(I)}(e_I(I))=e_I(I)\cup (e_I(I))^*$,
so that the map $e_{(A,I)}$ is defined by the following formula:
for every $a\in I$, $e_{(A,I)}(a^*)=(e_I(a))^*$.) Set
$I\ap=e_I(I)$ and $A\ap=e_{(A,I)}(A)$. Then the map
$e_{(A,I)}:(A,I)\lra(A\ap,I\ap)$ is a $\PLBA$-isomorphism. Note
that $(A\ap,I\ap)=(E^z\circ E^p)(A,I)$. Hence,
$e_{(A,I)}:(A,I)\lra(E^z\circ E^p)(A,I)$ is a $\PLBA$-isomorphism.
We will show that $e: Id_{\PLBA}\lra E^z\circ E^p$, defined by
$e(A,I)=e_{(A,I)}$ for every $(A,I)\in\card{\PLBA}$, is the
required natural isomorphism. Indeed, if $\p\in\PLBA((A,I),(B,J))$
and $\p\ap=(E^z\circ E^p)(\p)$ then we have to prove that
$e_{(B,J)}\circ \p=\p\ap\circ e_{(A,I)}$. Clearly, for doing this
it is enough to show that
$e_J\circ(\p_{|I})=(\p\ap)_{|e_I(I)}\circ e_I$. Since this is
obvious, we obtain that the functors $Id_{\PLBA}$ and $E^z\circ
E^p$ are naturally isomorphic.

Let us proceed with the composition $E^p\circ E^z$. Let $(B,I)$ be
a ZLBA. Then, by Corollary \ref{stfact}, the map
$\SI_{(B,I)}:Si(I)\lra B$, where $\SI_{(B,I)}(J)=\bv_B \{a\st a\in
J\}$ for every $J\in Si(I)$, is
 a Boolean isomorphism. We will show that $s: Id_{\PZLBA}\lra E^p\circ
E^z$, defined by $s(B,I)=(\SI_{(B,I)})\inv$ for every
$(B,I)\in\card{\PZLBA}$, is the required natural isomorphism.
Indeed, if $\p\in\PZLBA((A,I),(B,J))$ and $\p\ap=(E^p\circ
E^z)(\p)$ then we have to prove that $\SI_{(B,J)}\circ
\p\ap=\p\circ \SI_{(A,I)}$. Let $I_1\in Si(I)$. Then $(\p\circ
\SI_{(A,I)})(I_1)=\p(\bv_A I_1)$ and $(\SI_{(B,J)}\circ
\p\ap)(I_1)=\SI_{(B,J)}(\p\ap(I_1))=\SI_{(B,J)}(\bv_{Si(J)}\{\downarrow(\p(a))\st
a\in I_1\})=\bv_B\{\SI_{(B,J)}(\downarrow(\p(a)))\st a\in I_1\}=
\bv_B\p(I_1)$.   So, we have to prove that $\p(\bv_A I_1)=\bv_B
\p(I_1)$. Set $b=\p(\bv_A I_1)$ and $c=\bv_B \p(I_1)$. Since
$a\le\bv_A I_1$, for every $a\in I_1$, we have that $\p(a)\le b$
for every $a\in I_1$. Hence $c\le b$. We will now prove  that
$b\le c$. Since $J$ is dense in $B$, we get that $b=\bv_B\{d\in
J\st d\le b\}$. By condition (LBA), for every $d\in J$ there
exists $e_d\in I$ such that $d\le\p(e_d)$. So, let $d\in J$ and
$d\le b$. Since $I_1\vee\neg I_1=I$, there exist $e_d^1\in I_1$
and $e_d^2\in\neg I_1$ such that $e_d=e_d^1\vee e_d^2$. Now we
obtain that $d\le \p(e_d)\we b=\p(e_d\we\bv_A
I_1)=\p(\bv_A\{e_d\we a\st a\in I_1\})=\p(\bv_A\{e_d^1\we a\st
a\in I_1\})=\p(e_d^1\we\bv_A I_1)\le\p(e_d^1)\le c$. Thus
$b=\bv_B\{d\in J\st d\le b\}\le c$. So, the functors $Id_{\PZLBA}$
and $E^p\circ E^z$ are naturally isomorphic. \sqs

\begin{cor}\label{zlcbij}
 There exists a bijective correspondence between the
classes of all (up to $\PLBA$-isomorphism) PLBAs, all (up to
$\ZLBA$-isomorphism) ZLBAs  and all (up to homeomorphism) locally
compact zero-dimensional Hausdorff spaces.
\end{cor}

We can even express Theorem \ref{genstone} in a more simple form
which is very close to the results obtained by M. Stone in
\cite{ST}:

\begin{theorem}\label{genstonegst}
The category\/ $\PZLC$  is dually equivalent to the category\/
$\GBPL$ whose objects are all generalized Boolean pseudolattices
and whose morphisms are all $\{0\}$-pseudolattice homomorphisms
between its objects satisfying condition (LBA).
\end{theorem}

\doc By virtue of Theorem \ref{genstone}, it is enough to show
that the categories $\GBPL$ and $\PLBA$ are equivalent.

Define a functor $E^l:\PLBA\lra\GBPL$ by setting $E^l(A,I)=I$, for
every $(A,I)\in\card{\PLBA}$, and for every
$\p\in\PLBA((A,I),(B,J))$, put $E^l(\p)=\p_{|I}:I\lra J$. Using
Fact \ref{gbplplba}(a) and condition (PLBA), we get that $E^l$ is
a well-defined functor.

Define a functor $E^g:\GBPL\lra \PLBA$ by setting
$$E^g(I)=(B_{Si(I)}(e_I(I)),e_I(I))$$ for every $I\in\card{\GBPL}$
(see \ref{gbpapi}(b) and \ref{day} for the notations), and for
every $\p\in\GBPL(I,J)$ define $E^g(\p):B_{Si(I)}(e_I(I))\lra
B_{Si(J)}(e_J(J))$ to be the obvious extension of the map
$\p_e:e_I(I)\lra e_J(J)$ defined by $$\p_e(\dar(a))=\dar(\p(a)).$$
Then, using Facts \ref{gbpapi}(a) and \ref{gbplplba}(b), it is
easy to see that $E^g$ is a well-defined functor.

Finally, it is almost obvious that the compositions $E^g\circ E^l$
and $E^l\circ E^g$ are naturally isomorphic to the corresponding
identity functors. \sqs

\begin{cor}\label{zlcbijs}{\rm (M. Stone \cite{ST})}
 There exists a bijective correspondence between the
class of all (up to $\GBPL$-isomorphism) generalized Boolean
pseudolattices and all (up to homeomorphism) locally compact
zero-dimensional Hausdorff spaces.
\end{cor}

Note that in \cite{ST1}, M. Stone proves that there exists a
bijective correspondence between generalized Boolean
pseudolattices and Boolean rings (with or without unit).



\section{A description of  $\DLC$-products of LCAs}

\begin{defi}\label{prodclca}
\rm Let $\GA$ be a set and
$\{(A_\g,\rho_\g,\BBBB_\g)\st\g\in\GA\}$ be a family of LCAs. Let
$A=\prod\{A_\g\st\g\in\GA\}$ be the product of the  Boolean
algebras $\{A_\g\st\g\in\GA\}$ in the category $\Bo$ of Boolean
algebras and Boolean homomorphisms (i.e., $A$ is the Cartesian
product of the family $\{A_\g\st\g\in\GA\}$, construed as a
Boolean algebra with respect to the coordinate-wise operations).
Let $\BBBB=\{(b_\g)_{\g\in\GA}\in\prod\{\BBBB_\g\st\g\in\GA\}\st
\card{\{\g\in \GA\st b_\g\neq 0\}}<\aleph_0\}$, where
$\prod\{\BBBB_\g\st\g\in\GA\}$
 is the Cartesian product of the family $\{\BBBB_\g\st\g\in\GA\}$
 (in other words,
 $\BBBB$ is the $\s$-product of the family
$\{\BBBB_\g\st\g\in\GA\}$ with base point $0=(0_\g)_{\g\in\GA}$).
For any two points $a=(a_\g)_{\g\in\GA}\in A$ and
$b=(b_\g)_{\g\in\GA}\in A$, set $a\rho b$ if there exists
$\g\in\GA$ such that $a_\g\rho_\g b_\g$. Then the triple
$(A,\rho,\BBBB)$ is called a {\em product of the family of LCAs}
$\{(A_\g,\rho_\g,\BBBB_\g)\st\g\in\GA\}$. We will write
$(A,\rho,\BBBB)=\prod\{(A_\g,\rho_\g,\BBBB_\g)\st\g\in\GA\}$.
\end{defi}

\begin{fact}\label{prodf}
The product $(A,\rho,\BBBB)$ of a family
$\{(A_\g,\rho_\g,\BBBB_\g)\st\g\in\GA\}$
 of LCAs
 is an LCA.
\end{fact}

\doc The proof is straightforward. \sqs

\begin{pro}\label{prodclcacat}
 Let $\GA$ be a set and
$\{(A_\g,\rho_\g,\BBBB_\g)\st\g\in\GA\}$ be a family of CLCAs.
Then the source
$\{\pi_\g:(A,\rho,\BBBB)\lra(A_\g,\rho_\g,\BBBB_\g)\st\g\in\GA\}$,
where $(A,\rho,\BBBB)=\prod\{(A_\g,\rho_\g,\BBBB_\g)\st\g\in\GA\}$
(see \ref{prodclca}) and,  for every $a=(a_\g)_{\g\in\GA}\in A$
and every $\g\in\GA$, $\pi_\g(a)=a_\g$, is a product of the family
$\{(A_\g,\rho_\g,\BBBB_\g)\st\g\in\GA\}$ in the category $\DLC$.
\end{pro}

\doc By Fact \ref{prodf}, $(A,\rho,\BBBB)$ is an LCA and since $A$
is a complete Boolean algebra, we get that  $(A,\rho,\BBBB)$ is a
CLCA. It is easy to see that, for every $\g\in\GA$, $\pi_\g$ is a
$\DLC$-morphism.

Let $X_\g=\LAM^a(A_\g,\rho_\g,\BBBB_\g)$ for every $\g\in\GA$, and
let $X=\bigoplus\{X_\g\st \g\in\GA\}$ be the topological sum of
the family $\{X_\g\st \g\in\GA\}$. Then the sink of inclusions
$\{i_\g:X_\g\lra X\st\g\in\GA\}$ is a coproduct in the category
$\HLC$ (briefly, $\HLC$-coproduct) of the family $\{X_\g\st
\g\in\GA\}$. Since $\LAM^t$  is a duality (by  \cite[Theorem
2.14]{D5}), the source
$\PP=\{\LAM^t(i_\g):\LAM^t(X)\lra\LAM^t(X_\g)\st\g\in\GA\}$ is a
$\DLC$-product of the family $\{\LAM^t(X_\g)\st \g\in\GA\}$. Then,
clearly, the source
$\QQ=\{(\l^g_{A_\g})\inv\di\LAM^t(i_\g):\LAM^t(X)\lra
(A_\g,\rho_g,\BBBB_\g)\st\g\in\GA\}$ is a $\DLC$-product of the
family $\{(A_\g,\rho_\g,\BBBB_\g)\st\g\in\GA\}$. Set
$\a_\g=(\l^g_{A_\g})\inv\di\LAM^t(i_\g)$. We will show that there
exists a $\DLC$-isomorphism $\a:\LAM^t(X)\lra (A,\rho,\BBBB)$ such
that, for any $\g\in\GA$, $\pi_\g\di\a=\a_\g$. Obviously, this
will imply that the source
$\{\pi_\g:(A,\rho,\BBBB)\lra(A_\g,\rho_g,\BBBB_\g)\st\g\in\GA\}$
is a $\DLC$-product of the family
$\{(A_\g,\rho_\g,\BBBB_\g)\st\g\in\GA\}$. Set, for every $F\in
RC(X)$ and any $\g\in\GA$, $F_\g=F\cap X_\g$. Then $F_\g\in
RC(X_\g)$ for every  $\g\in\GA$. Define the map $\a:RC(X)\lra A$
by $\a(F)=((\l^g_{A_\g})\inv(F_\g))_{\g\in\GA}$, for every $F\in
RC(X)$. Since $\LAM^t(X)=RC(X)$ and $\LAM^t(X_\g)=RC(X_\g)$, it is
easy to see that the map $\a$ is a $\DLC$-isomorphism between
$\LAM^t(X)$ and $(A,\rho,\BBBB)$. Further, for any $\g\in\GA$ and
any $F\in RC(X)$,
$\LAM^t(i_\g)(F)=\cl_{X_\g}(i_\g\inv(\int_X(F)))$ (see
\cite[Theorem 2.14]{D5}). We get that $\LAM^t(i_\g)(F)=F_\g$ which
implies easily that $\pi_\g\circ\a=\a_\g$, for every $\g\in\GA$.
Thus, by (DLC5), $\pi_\g\di\a=\a_\g$, for every $\g\in\GA$. \sqs

\section{The notion of weight  of an LCA}

The next definition and proposition generalize the analogous
definition and statement of de Vries \cite{dV}. Note that our
$``$\/\/base" (see the definition below) appears in \cite{dV} (for
NCAs) as $``$dense set". (See \cite[Definition 1.1]{D5} for the
notion $``$NCA").

\begin{defi}\label{clcawe}
\rm Let $(A,\rho,\BBBB)$ be an LCA and $\BB$ be a subset of
$\BBBB$. Then $\BB$ is called a {\em base} (or a {\em dV-dense
subset}) of $(A,\rho,\BBBB)$ if for each $a,c\in\BBBB$ such that
$a\llx c$ there exists $b\in\BB$ with $a\le b\le c$. The cardinal
number $w(A,\rho,\BBBB)=\min\{\card\BB\st \BB$ is a base of
$(A,\rho,\BBBB)\}$ is called a {\em weight of} $(A,\rho,\BBBB)$.
\end{defi}

\begin{fact}\label{wfact}
 If $(A,\rho,\BBBB)$ is an LCA and $\BB$ is a subset of
$\BBBB$ then $\BB$ is a  base of $(A,\rho,\BBBB)$ iff for each
$a,c\in\BBBB$ such that $a\llx c$ there exists $b\in\BB$ with
$a\llx b\llx c$.
\end{fact}

\doc ($\Rightarrow$) Let $a,c\in\BBBB$ and $a\llx c$. Then, by (BC1), there exists $d,e\in\BBBB$ with
$a\llx d\llx e\llx c$. Now, there exists $b\in\BB$ such that $d\le
b\le e$. Therefore $a\llx b\llx c$.

\noindent($\Leftarrow$) This is clear. \sqs

\begin{pro}\label{weightpiw}
Let $\tau$ be an infinite cardinal number, $(A,\rho,\BBBB)$ be an
LCA and $X=\Psi^a(A,\rho,\BBBB)$. Then $w(X)=\tau$ iff
$w(A,\rho,\BBBB)=\tau$.
\end{pro}

\doc We know that the family $\BB_0=\{\int_X(\lag(a))\st
a\in\BBBB\}$ is a base of $X$.

   Let $w(X)=\tau$. Then there exists a base $\BB\ap$ of $X$ such that
$\BB\ap\sbe\BB_0$ and $\card{\BB\ap}=\tau$. Let $\BB$ be the
sub-join-pseudolattice of $\BBBB$  generated  by the set
$\{a\in\BBBB\st\int(\lag(a))\in\BB\ap\}$. It is clear that $\BB$
is a base of $(A,\rho,\BBBB)$. Hence, $w(X)\ge w(A,\rho,\BBBB)$.

Conversely, if $\BB$ is a base of $(A,\rho,\BBBB)$ and
$\card{\BB}=\tau$, then it is easy to see that
$\BB\ap=\{\int(\lag(a))\st a\in\BB\}$ is a base of $X$. Thus,
$w(X)\le w(A,\rho,\BBBB)$. \sqs

\begin{cor}\label{auth}
Let $\BB$ be a base of an LCA $(A,\rho,\BBBB)$ with infinite
weight. Then there exists a base $\BB_1$ of $(A,\rho,\BBBB)$ such
that $\BB_1\sbe\BB$ and $\card{\BB_1}=w(A,\rho,\BBBB)$.
\end{cor}

\doc This follows from the second part of the proof of Theorem
\ref{weightpiw} and the well-known Alexandroff-Urysohn Theorem for
bases (see, e.g., \cite[Theorem 1.1.15]{E}). \sqs

\begin{theorem}\label{metrclca}
Let $X\in\card\HLC$. Then $X$ is metrizable iff there exists a
set\/ $\GA$ and a family $\{(A_\g,\rho_\g,\BBBB_\g)\st\g\in\GA\}$
of CLCAs such that
$\LAM^t(X)=\prod\{(A_\g,\rho_\g,\BBBB_\g)\st\g\in\GA\}$ and, for
each $\g\in\GA$, $w(A_\g,\rho_\g,\BBBB_\g)\le\aleph_0$.
\end{theorem}

\doc The following theorem is well-known (see, e.g., \cite{A} or
the more general theorem \cite[Theorem 5.1.27]{E}): a locally
compact Hausdorff space is metrizable iff it is a topological sum
of locally compact Hausdorff spaces with countable weight. Since,
by  \cite[Theorem 2.14]{D5}, $\LAM^t$ is a duality functor, it
converts the $\HLC$-sums in $\DLC$-products. Hence, our assertion
follows from the cited above theorem and Propositions
\ref{prodclcacat} and \ref{weightpiw}. \sqs

\begin{nota}\label{as}
\rm  Let $(A,\rho,\BBBB)$ be an LCA. We set
$$(A,\rho,\BBBB)_S=\{a\in A\st a\ll_\rho a\}.$$ We will  write
simply $``A_S$" instead of $``(A,\rho,\BBBB)_S$" when this does
not leads to an ambiguity.
\end{nota}

\begin{pro}\label{nuldim}
Let  $(A,\rho,\BBBB)$ be an LCA. Then the space
$\Psi^a(A,\rho,\BBBB)$ is zero-dimensional iff the set
$A_S\cap\BBBB$ is a base of $(A,\rho,\BBBB)$.
\end{pro}

\doc  Let $X=\Psi^a(A,\rho,\BBBB)$. Then the family $\{\int_X(\lag(a))\st
a\in\BBBB\}$ is a base of $X$.

\noindent($\Rightarrow$)  If $U$ is a clopen compact subset of $X$
then, clearly, $U=\lag(a)$ for some $a\in\BBBB\cap A_S$. This
implies that $A_S\cap\BBBB$ is a base of $(A,\rho,\BBBB)$.

\noindent($\Leftarrow$) Let $x\in X$ and $U$ be a neighborhood of
$x$. Then there exist $a,b\in\BBBB$ such that
$x\in\int(\lag(a))\sbe\lag(a)\sbe\int(\lag(b))\sbe U$. Then $a\ll
b$. Hence there exists $c\in A_S\cap\BBBB$ such that $a\le c\le
b$. Then $V=\lag(c)$ is clopen in $X$ and $x\in V\sbe U$. \sqs

In the sequel, we will denote by $C$  the Cantor set.

 Note that
$RC(C)$ is isomorphic to the minimal completion $B$ of a free
Boolean algebra $A$ with $\aleph_0$ generators or, equivalently,
$RC(C)$ is the unique (up to isomorphism)  atomless complete
Boolean algebra $B$  containing a countable dense subalgebra $A$
(see, e.g., \cite{Dw}).   Defining in $B$ a relation $\rho$ by
$a(-\rho)b$ (where $a,b\in B$) iff there exists  $c\in A$ such
that $a\le c\le b^*$, we get that $(B,\rho)$ is a CNCA
CA-isomorphic to the CNCA $(RC(C),\rho_C)$ (see \cite[Definition
1.1]{D5} for these notions). We will now obtain a generalization
of this construction.

\begin{pro}\label{nuldimlm}
Let $A_0$ be a dense  Boolean subalgebra of a Boolean algebra $A$.
Then setting, for every $a,b\in A$,  $a\ll_\rho b$ if there exists
 $c\in A_0$ such that $a\le c\le b$, we obtain
a normal contact relation $\rho$ on $A$ such that
$(A,\rho)_S=A_0$, $A_0$ is the smallest base of $(A,\rho)$ and
$w(A,\rho)=\card{A_0}$. Also, $\Psi^a(A_0,\rho_s)=S^a(A_0)$ (where
$S^a:\Bool\lra\ZHC$ is the Stone duality functor) and when $A$ is
complete then $S^a(A_0)=\Psi^a(A,\rho)$.
\end{pro}

\doc  It is easy to check that the
relation $\rho$ satisfies conditions ($\ll 1$)-($\ll 7$). (For
establishing ($\ll 5$) and ($\ll 6$) use the fact that for every
$c\in A_0$ we have, by the definition of the relation $\llx$, that
$c\llx c$.) The rest is clear. (Note only that if $A$ is complete
and $X=S^a(A_0)$ then the NCAs $(A,\rho)$ and $(RC(X),\rho_X)$ are
NCA-isomorphic.) \sqs

\section{On the  $\pi$-weight  of a poset}

\begin{defi}\label{pibaseposet}
\rm  Let $(A ,\le)$ be a poset. We set $\pi
w(A,\le)=\min\{\card\BB\st \BB$ is dense in $(A,\le)\}$; the
cardinal number $\pi w(A,\le)$ is called a $\pi$-{\em weight of
the poset} $(A,\le)$.
\end{defi}

The term {\em density of} $(A,\le)$ instead that of $\pi$-{\em
weight} is usually used. Our reason for introducing a new term is
Proposition \ref{piweightpisem} which is proved below.

Obviously, (BC3) and (BC1) imply that every base of a local
contact algebra $(A,\rho,\BBBB)$ is a dense subset of $A$. Hence,
for every LCA $(A,\rho,\BBBB)$, $\pi w(A)\le w(A,\rho,\BBBB)$.

\begin{fact}\label{piwel}
Let $(A,\rho,\BBBB)$ be an LCA and $\BB$ be a subset of $A$. Then
the following conditions are equivalent:

\noindent(a) $\BB$ is a dense subset of $(A,\rho,\BBBB)$;

\noindent(b) for each $a\in A\stm\{0\}$ there exists
$b\in\BB\stm\{0\}$ such that $b\llx a$;

\noindent(c) for each $a\in\BBBB\stm\{0\}$, $a=\bv\{b\in\BB\st
b\llx a\}$;

\noindent(d) for each $a\in A\stm\{0\}$, $a=\bv\{b\in\BB\st b\llx
a\}$.
\end{fact}

\doc The implications (a) $\leftrightarrow$ (b), (c) $\leftrightarrow$ (d)
and (d) $\rightarrow$ (a) are clear. We need only to show that (a)
$\rightarrow$ (d).

Let $a\in A\stm\{0\}$. Then $a=\bv\{b\in\BB\st b\le a\}$. Let
$a_1\in A$ and $a_1\ge b$ for every $b\in\BB$ such that $b\llx a$.
Suppose that $a_1\not\ge a$. Then $a\we a_1^*\neq 0$. By (BC3),
there exists $c\in\BBBB\stm\{0\}$ such that $c\llx a\we a_1^*$.
There exists $b\in\BB\stm\{0\}$ with $b\le c$. Then $b\llx a\we
a_1^*$. Thus $b\llx a$ and hence $b\le a_1$. Therefore $b\le
a_1\we a_1^*=0$, a contradiction. So, $a=\bv\{b\in\BB\st b\llx
a\}$. \sqs

The next fact is obvious.

\begin{fact}\label{piwelf}
Let $(A,\rho,\BBBB)$ be an LCA and $\BB$ be a dense subset of $A$.
Then $\BB\cap\BBBB$ is a dense subset of $A$.
\end{fact}

Recall that if $(X,\TT)$ is a topological space then: a) a family
$\BB$ of open subsets of  $(X,\TT)$ is called a $\pi$-{\em base}\/
of $(X,\TT)$ if for each $U\in\TT\stm\{\ems\}$ there exists $V\in
\BB\stm\{\ems\}$ such that $V\sbe U$; b) the cardinal number $\pi
w(X)=\min\{\card{\BB}\st \BB$ is a $\pi$-base of $(X,\TT)\}$ is
called a $\pi$-{\em weight}\/ of $(X,\TT)$.

\begin{defi}\label{pisemsp}
\rm A topological space $(X,\TT)$ is called  $\pi$-{\em
semiregular} if  the family $RO(X)$ is a $\pi$-base of $X$.
\end{defi}

Clearly, every semiregular space is $\pi$-semiregular. The
converse is not true. Indeed, it is easy to see that the space $X$
from \cite[Example 78]{SS} (known as $``$\/half-disc topology") is
a $\pi$-semiregular $T_{2\frac{1}{2}}$-space which is not
semiregular. On the other hand, if $X$ is an infinite set with the
cofinite topology then $X$ is not a $\pi$-semiregular space since
$RO(X)=\{\ems, X\}$.

We will now need a simple lemma.

\begin{lm}\label{dislm}
If\/ $\BB$ is a $\pi$-base of a space $X$ then there exists a
$\pi$-base $\BB\ap$ of $X$ such that $\BB\ap\sbe\BB$ and
$\card{\BB\ap}=\pi w(X)$.
\end{lm}

\doc Let $\BB_0$ be a $\pi$-base of $X$ with $\card{\BB_0}=\pi
w(X)$. Then for every non-empty $U\in\BB_0$ there exists
$V_U\in\BB\stm\{\ems\}$ such that $V_U\sbe U$. Obviously,
$\BB\ap=\{V_U\st U\in\BB_0\}$ is the required $\pi$-base. \sqs

\begin{pro}\label{piweightpisem}
If $X$ is a $\pi$-semiregular topological space, then  $\pi w(X)
=\pi w(RC(X))$.
\end{pro}

\doc Since $X$ is $\pi$-semiregular, $RO(X)$ is a $\pi$-base of $X$.
Hence, by Lemma \ref{dislm},   there exists a $\pi$-base $\BB$ of
$X$ such that $\BB\sbe RO(X)$ and $\card{\BB}=\pi w(X)$.
Obviously, $\BB$ is a dense subset of $(RO(X),\sbe)$ as well.
Hence,  $\pi w(X)\ge\pi w(RO(X))$. Clearly,  $\pi w(X) \le\pi
w(RO(X))$. Finally, note that $(RO(X),\sbe)$ and $(RC(X),\sbe)$
are isomorphic posets. \sqs

The  assertion which follows should be known. We will use it for
obtaining some slight generalizations of two results of V. I.
Ponomarev \cite{P1}.

\begin{lm}\label{ponlm}
Let $A$ and $B$ be Boolean algebras and $\p:A\lra B$ be a
function.

\noindent a) If $\p$  satisfies the following conditions:

\medskip

\noindent 1) $\p(a\vee b)=\p(a)\vee\p(b)$, for all $a,b\in A$,

\noindent 2) $\p(0_A)=0_B$ and $\p\inv(1_B)=1_A$,

\noindent 3) $\p(A)$ is dense in $B$,

\medskip

\noindent then the map $\p$ is a Boolean embedding (= injective
Boolean homomorphism) and $\pi w(A)=\pi w(B)$;

\medskip

\noindent b) If $A$ is complete then $\p:A\lra B$ is a Boolean
isomorphism iff $\p$ satisfies conditions 1)-3) from a).
\end{lm}

\doc a) Note that for every $a\in A$, $\p(a^*)\ge (\p(a))^*$.
Indeed, this follows from the equations $1_B=\p(1_A)=\p(a\vee
a^*)=\p(a)\vee\p(a^*)$. Further, let $a,b\in A$ and $\p(a)=\p(b)$.
Then $\p(a)\vee(\p(b))^*=1$. Hence $\p(a)\vee\p(b^*)=1$. Thus
$a\vee b^*=1$, i.e. $b\le a$. Analogously, starting with
$\p(b)\vee(\p(a))^*=1$, we get that $a\le b$. So, $\p$ is an
injection.

Let $\p(a)\le\p(b)$. We will show that then $a\le b$. Indeed, we
have that  $\p(a)\vee\p(b)=\p(b)$. Hence $\p(a\vee b)=\p(b)$. Thus
$a\vee b=b$, i.e. $a\le b$. Further, if $\p(a)\we\p(a^*)\neq 0$
then, by the density of $\p(A)$ in $B$, there exists $b\in A$ such
that $0\neq\p(b)\le\p(a)\we\p(a^*)$. Then $0\neq b\le a$ and $b\le
a^*$, i.e. $b=0$, a contradiction. Thus $\p(a)\we\p(a^*)=0$. Hence
$\p(a^*)\le(\p(a))^*$, which implies that $\p(a^*)=(\p(a))^*$. So,
$\p$ is a Boolean embedding. Since $\p(A)$ is dense in $B$, it is
easy to see that $\pi w(\p(A))\ge\pi w(B)$. Conversely, if $\BB$
is a dense subset of $B$ with $\card{\BB}=\pi w(B)$ then for every
$b\in \BB\stm\{0\}$ there exists $a_b\in A$ such that
$0\neq\p(a_b)\le b$. Set $\BB\ap=\{\p(a_b)\st b\in\BB\}$. Then,
clearly, $\BB\ap$ is dense in $\p(A)$. Therefore, $\pi w(A)=\pi
w(B)$.

b) By a), we need only to show that $\p$ is a surjection. This is
so because $A$ is complete and $\p(A)$ is dense in $B$ (see
\cite{Si}). \sqs

Recall the Ponomarev's result \cite{P} that a map $f:(X,\TT)\lra
(Y,\OO)$ is closed and irreducible iff it is a surjection and, for
every $U\in\TT\stm\{\ems\}$, $f^\sharp(U)\in\OO\stm\{\ems\}$.

\begin{defi}\label{quasipi}
\rm Let $f:(X,\TT)\lra (Y,\OO)$ be a continuous map. We will say
that $f$ is a $\pi$-{\em map}\/ if it is a closed irreducible map.
The map $f$ is called a {\em quasi-$\pi$-map} (respectively, an
{\em  MR-map})\/ if $\cl(f(X))=Y$ and for every
$U\in\TT\stm\{\ems\}$ (respectively, for every $U\in
RO(X)\stm\{\ems\}$) we have that $\int(f^\sharp(U))\nes$.
\end{defi}

The name $``$quasi-$\pi$-map" is chosen because the definition of
these maps is similar to the definition of  quasi-open maps.  As
we shall see later, our MR-maps almost coincide with the
continuous irreducible in the sense of Mioduszewski and Rudolf
\cite{MR} maps.

 Obviously, every $\pi$-map is a quasi-$\pi$-map and every
 quasi-$\pi$-map is an MR-map.  If $X$ is $\pi$-semiregular then
 every MR-map $f:X\lra Y$ is a
 quasi-$\pi$-map.
 Since, clearly, the
 dense embeddings are quasi-$\pi$-maps, we get that not every quasi-$\pi$-map is a
 $\pi$-map. It can be easily shown that the composition of two
 quasi-$\pi$-maps is a quasi-$\pi$-map.

\begin{fact}\label{quasipif}
A continuous map $f:X\lra Y$ is a quasi-$\pi$-map (respectively,
an MR-map) iff $\cl(f(X))=Y$ and $cl(f(F))\neq Y$ for each closed
proper subset $F$ of $X$ (respectively, for each $F\in
RC(X)\stm\{X\}$).
\end{fact}

\doc It is well-known that for every subset $M$ of $X$,
$f^\sharp(M)=Y\stm f(X\stm M)$; hence $\int(f^\sharp(M))=Y\stm
\cl(f(X\stm M))$. The rest is clear. \sqs

\begin{cor}\label{quasipicl}
A  closed map is a quasi-$\pi$-map iff it is a $\pi$-map.
\end{cor}

A surjective map $f:X\lra Y$, where $Y$ is a Hausdorff space, is
{\em irreducible in the sense of} Mioduszewski and  Rudolf
\cite{MR} if for every $F\in RC(X)$, $F\neq X$ implies that
$\cl(f(F))\neq Y$. Hence, the only difference between MR-maps and
continuous irreducible maps in the sense of \cite{MR} is that
MR-maps are not assumed to be surjections and $Y$ is not assumed
to be Hausdorff. As it is noted in \cite{MR}, if $X$ is compact
then every irreducible in the sense of \cite{MR} continuous map
$f:X\lra Y$ is an irreducible map.

\begin{lm}\label{skellm1}
A continuous function $f:X\lra Y$ is  skeletal if and only if\/
$\int(\cl(f(U)))\nes$, for every  non-empty regular open subset
$U$ of $X$.
\end{lm}

\doc By \cite[Lemma 2.4]{D1}, a function $f:X\lra Y$ is  skeletal
if and only if\/ $\int(\cl(f(U)))\nes$, for every  non-empty  open
subset $U$ of $X$. Hence, we need only to show that if $f$ is
continuous and $\int(\cl(f(U)))\nes$ for every  non-empty regular
open subset $U$ of $X$, then $\int(\cl(f(U)))\nes$ for every
non-empty open subset $U$ of $X$. Let $U$ be an open non-empty
subset of $X$. Set $U\ap=\int(\cl(U))$. Then $U\ap\in RO(X)$,
$U\ap\nes$ and $\cl(U\ap)=\cl(U)$. Using continuity of $f$, we get
that $\cl(f(U))=\cl(f(\cl(U)))=\cl(f(\cl(U\ap)))=\cl(f(U\ap))$.
\sqs

\begin{pro}\label{skelpro}
Every MR-map is skeletal.
\end{pro}

\doc Let $f:X\lra Y$ and $U\in RO(X)\stm\{\ems\}$. We will show that
$\int(f^\sharp(U))\sbe\int(\cl(f(U)))$. Then Lemma \ref{skellm1}
will imply that $f$ is a skeletal map. Let
$y\in\int(f^\sharp(U))$. Then there exists an open neighborhood
$O$ of $y$ such that $f\inv(O)\sbe U$. Then $O\cap
f(X)=f(f\inv(O))\sbe f(U)$. Since $f(X)$ is dense in $Y$, we get
that $O\sbe \cl(O)=\cl(O\cap f(X))\sbe \cl(f(U))$. Hence,
$y\in\int(\cl(f(U)))$. \sqs

\begin{pro}\label{rciso1}
Let $f:X\lra Y$ be an MR-map. Then the Boolean algebras $RC(X)$
and $RC(Y)$ are isomorphic. If, moreover, $X$ and $Y$ are
$\pi$-semiregular spaces, then $\pi w(X)=\pi w(Y)$.
\end{pro}

\doc Since the map $f$ is skeletal, we have, by \cite[Lemma 2.6]{D1},
that for every $F\in RC(X)$, $\cl(f(F))\in RC(Y)$. Now, define a
map $\p:RC(X)\lra RC(Y)$ by $\p(F)=\cl(f(F))$, for every $F\in
RC(X)$. Obviously, $\p$ satisfies conditions 1) and 2) of Lemma
\ref{ponlm} (see Fact \ref{quasipif}). Further, let $G\in RC(Y)$
and $G\nes$. Set $F=\cl(f\inv(\int(G)))$. Then, clearly, $F\in
RC(X)$ and since $\cl(f(X))=Y$, we have that $F\nes$. The
continuity of $f$ implies that $f(F)\sbe G$. Thus $\p(F)\sbe G$.
Hence $\p(RC(X))$ is dense in $RC(Y)$. Therefore we get, by Lemma
\ref{ponlm}(b), that the Boolean algebras $RC(X)$ and  $RC(Y)$ are
isomorphic. Now, Proposition \ref{piweightpisem} implies that $\pi
w(X)=\pi w(Y)$. \sqs

\begin{cor}\label{corpon1}{\rm (\cite{P})}
 If $f:X\lra Y$ is a $\pi$-map then $RC(X)$ and $RC(Y)$ are
isomorphic Boolean algebras.
\end{cor}

 Obviously, Proposition \ref{rciso1} implies also
 (in the class of $\pi$-semiregular spaces) the result of
Ponomarev \cite{P} that if $Y$ is an image of $X$ under a
$\pi$-map then $\pi w(X)= \pi w(Y)$.

\section{Co-absolute spaces}

\begin{pro}\label{prorho}
Let $A$ be a Boolean algebra and  $\pi w(A)\ge\aleph_0$. Then
there exists a normal contact relation $\rho$ on $A$ such that
$w(A,\rho)=\pi w(A)$ and $(A,\rho)_S$ is a base of $(A,\rho)$.
\end{pro}

\doc  Let $\pi w(A)=\tau$.  Then there exists a dense subset $B_0$ of $A$
with $\aleph_0\le\card{B_0}=\tau$. Let $B$ be the Boolean
subalgebra of $A$ generated by $B_0$. Now, Proposition
\ref{nuldimlm} implies that there exists a normal contact relation
$\rho$ on $A$ such that $(A,\rho)_S$ is a base of $(A,\rho)$ and
$w(A,\rho)=\card{B}$. Since $\card{B} =\card{B_0}=\tau$, we get
that $w(A,\rho)=\tau$. \sqs

\begin{pro}\label{pon2}
Let $X$ be a $\pi$-semiregular space and $\pi w(X)\ge\aleph_0$.
Then  there exists a compact Hausdorff zero-dimensional space $Y$
with $w(Y)=\pi w(X)$ for which the Boolean algebras $RC(X)$ and
$RC(Y)$ are isomorphic.
\end{pro}

\doc   Let $\pi
w(X)=\tau$. Set $A=RC(X)$. Then, by \ref{piweightpisem}, $\pi
w(A)=\tau$. Hence, by Proposition \ref{prorho}, there exists a
normal contact relation $\rho$ on $A$ such that $w(A,\rho)=\tau$
and $(A,\rho)_S$ is a base of $(A,\rho)$. Thus, using Propositions
\ref{nuldim} and \ref{weightpiw}, we get that  $Y=\Psi^a(A,\rho)$
is a compact Hausdorff zero-dimensional space with $w(Y)=\tau$.
Finally, by de Vries Duality Theorem, $RC(Y)$ is isomorphic to
$A$, i.e. to $RC(X)$.
 \sqs

 We will give also a second proof of Proposition \ref{pon2}
 which uses only the Stone Duality Theorem and some well-known facts
 about  minimal completions: let
 $\pi w(X)=\tau$; then there exists a dense Boolean subalgebra $B$ of
 $RO(X)$ with $\card{B}=\tau$; further, $RO(X)$ is a minimal
 completion of $B$; let $Y=S^a(B)$; then $Y$ is a zero-dimensional
 compact Hausdorff space with $CO(Y)\cong B$; hence $w(Y)=\tau$;
 since $RO(Y)$ is a minimal completion of $CO(Y)$, we get that
 $RO(Y)\cong RO(X)$.

In connection with Proposition \ref{pon2}, let us mention a fact
which follows immediately from Stone Duality Theorem:

\begin{fact}\label{pon2f}
For every topological space $X$ there exists a compact Hausdorff
extremally disconnected space $Y$ with $RC(Y)$ isomorphic to
$RC(X)$.
\end{fact}

\doc Set $Y=S^a(RC(X))$. Then  $Y$ is an
extremally disconnected compact Hausdorff space and $RC(Y)\cong
RC(X)$. \sqs

Note that for every infinite set $X$ with the cofinite topology on
it, the space $Y$ from Fact \ref{pon2f} is an one-point space;
thus, in general, there is no such connection  between the
$\pi$-weight of $X$ and the weight of $Y$ as in Proposition
\ref{pon2}.

We will now show that Proposition \ref{pon2} implies  Ponomarev's
theorem \cite{P1} that a  compact Hausdorff space $X$ is
co-absolute with a compact metric space iff $\pi w(X)\le\aleph_0$.

Recall first that if $X$ is a regular space then a space $EX$ is
called an {\em absolute} of $X$ iff there exists  a perfect
irreducible map $\pi_X:EX\lra X$ and every perfect  irreducible
 preimage of $EX$ is homeomorphic to $EX$ (see, e.g., \cite{PS}).
 Two regular spaces are said to be {\em co-absolute}\/ if their
 absolutes are homeomorphic.
 It
is well-known that:  a) the absolute is unique up to
homeomorphism; b) a space $Y$ is an absolute  of a regular space
$X$ iff $Y$ is an extremally disconnected Tychonoff space
 for which there exists a perfect irreducible map
$\pi_X:Y\lra X$;  c) if $X$ is a compact Hausdorff space then
$EX=S^a(RC(X))$, where $S^a$ is the Stone contravariant functor.
Taking the above statement b) as a definition of the absolute of a
regular space, we will give some new proofs of the existence and
the uniqueness of absolutes of locally compact Hausdorff spaces
and we will describe the dual objects (i.e. the images under the
contravariant functor  $\Psi^t$ (see \cite[(5) in the proof of
Theorem 2.1]{D5}) of these absolutes. For doing this we will need
a lemma which is contained in the proof of \cite[Theorem 2.11]{D1}
but  is not formulated explicitly there.

\begin{lm}\label{lead}
Let $f:X\lra Y$ be a skeletal map. Then the map $\psi:RC(X)\lra
RC(Y)$, defined by $\psi(F)=\cl(f(F))$ for every $F\in RC(X)$, is
a left adjoint to the map $\p:RC(Y)\lra RC(X)$ defined by
$\p(G)=\cl(f\inv(\int(G)))$ for every $G\in RC(Y)$ (i.e. $\psi$ is
the unique order preserving map from $RC(X)$ to $RC(Y)$ such that
for every $F\in RC(X)$, $F\sbe\p(\psi(F))$, and for every $G\in
RC(Y)$, $\psi(\p(G))\sbe G$).
\end{lm}

\doc See the beginning of the proof of \cite[Theorem 2.11]{D1}.
\sqs

A new proof of the existence of an absolute of a locally compact
Hausdorff space is given in the next proposition, where the dual
object of this absolute is described as well.

\begin{pro}\label{absolute}
Let $(A,\rho,\BBBB)$ be a CLCA and $X=\Psi^a(A,\rho,\BBBB)$. Then
 the space
$\Psi^a(A,\rho_s,\BBBB)$ (see \cite[Example 1.2]{D5} for the
notation $\rho_s$) is an  absolute of $X$.
\end{pro}

\doc Let $Y=\Psi^a(A,\rho_s,\BBBB)$. Then,  as it is shown in \cite[Proposition 2.14(b)]{D},
$Y$ is a locally compact Hausdorff extremally disconnected space.

Define a map $i:(A,\rho,\BBBB)\lra (A,\rho_s,\BBBB)$, by $i(a)=a$
for every $a\in A$. Obviously, the left adjoint to $i$ is
$j=i\inv(=i)$. Then, by  \cite[Proposition 2.26]{D5} and
\cite[Theorem 2.15]{D1}, we get that the map $f=\LAM^a(i):Y\lra X$
is a perfect skeletal map. Since, clearly, $i$ is a
$\PAL$-morphism, we get that $f$ is a surjection (see
\cite[Theorem 2.11]{D}). Let $\p$ and $\psi$ be defined as in
\ref{lead}. Then $\p=\LAM^t(f)$, and  \cite[Lemma 3.18]{D5}
implies that $\p$ is a Boolean isomorphism (because
$\p=\l^g_{(A,\rho_s,\BBBB)}\circ
i\circ(\l^g_{(A,\rho,\BBBB)})\inv$). Hence its left adjoint $\psi$
is also an isomorphism. Hence $\psi\inv(Y)=X$ (i.e.,
$\psi\inv(1)=1$). This means that $f$ is an irreducible map.
Therefore $Y$ is an absolute of $X$. \sqs

A new proof (using only our methods) of the uniqueness (up to
homeomorphism) of the absolute of a locally compact Hausdorff
space is given in the following proposition:

\begin{pro}\label{absoluteun}
The absolute of a locally compact Hausdorff space $X$ is unique up
to  homeomorphism.
\end{pro}

\doc Let $Y$ be an absolute of $X$, i.e.  $Y$ is an extremally
disconnected Tychonoff space and there exists a perfect
irreducible map $f:Y\lra X$. Then $Y$ is a locally compact space
(as a perfect preimage of a locally compact space) and $f$ is a
$\pi$-map. Let $\Psi^t(X)=(A,\rho,\BBBB)$ and
$\Psi^t(Y)=(B,\eta,\BBBB\ap)$. Then, by \cite[Proposition
2.14(b)]{D}, $\eta=\rho_s$ (see \cite[Example 1.2]{D5} for the
notation $\rho_s$). Since $f$ is a quasi-open (and, hence,
skeletal) map, we get, by \cite[Theorem 2.15]{D1}, that
$\p=\LAM^t(f):(A,\rho,\BBBB)\lra (B,\rho_s,\BBBB\ap)$ is a
complete Boolean homomorphism such that: a) if $a\in\BBBB$ then
$\p(a)\in\BBBB\ap$, and b) if $b\in\BBBB\ap$ then
$\psi(b)\in\BBBB$, where $\psi$ is the left adjoint of $\p$.
Moreover, since $f$ is a surjection, $\p$ is an injection (see
\cite[Theorem 2.11]{D}). By Lemma \ref{lead}, for every $F\in
RC(Y)(=B)$, $\psi(F)=f(F)$. In order to show that $\p$ is a
surjection, we need only to prove that for every $F\in B$,
$\p(\psi(F))\sbe F$ (then (see Lemma \ref{lead}) we will have that
$\p(\psi(F))= F$). So, let $F\in B$. Since $\int(f\inv(G))\spe
f\inv(\int(G))$ for every $G\sbe X$ (because $f$ is continuous),
it is enough to prove that $\int(f\inv(f(F)))\sbe F$. Let
$x\in\int(f\inv(f(F)))$ Then there exists an open neighborhood
$Ox$ of $x$ such that $Ox\sbe f\inv(f(F))$. Then, for every open
neighborhood $Vx$ of $x$ such that $Vx\sbe Ox$, $\ems\neq
f^\sharp(Vx)\sbe f(Vx)\sbe f(Ox)\sbe f(F)$. Let $y\in
f^\sharp(Vx)$. Then $f\inv(y)\sbe Vx$ and $y\in f(F)$. Hence there
exists $z\in F$ such that $f(z)=y$. Thus, $z\in f\inv(y)\sbe Vx$,
i.e. $z\in Vx\cap F$. Therefore, $x\in F$. So, $\p$ is a
bijection. Then $\psi=\p\inv$ and we get that
$\p(\BBBB)=\BBBB\ap$. Hence $\p:A\lra B$ is a Boolean isomorphism,
$\p(\BBBB)=\BBBB\ap$ and $\eta=\rho_s$. If $EX$ is the absolute of
$X$ constructed in Proposition \ref{absolute}, we get, by Roeper
Theorem (see  \cite[Theorem 2.1]{D5}),  that $Y$ is homeomorphic
to $EX$. \sqs

Now, our methods permit to obtain easily a slightly different form
of a well-known theorem of Ponomarev \cite{P1}.

\begin{theorem}\label{pon21}
Let $X$ be a compact Hausdorff space and $\tau$ be an infinite
cardinal number. If $\pi w(X)=\tau$ then $X$ is co-absolute with a
compact Hausdorff zero-dimensional space $Y$ with $w(Y)=\tau$.
\end{theorem}

\doc  By
Proposition \ref{pon2}, there exists a zero-dimensional compact
Hausdorff space $Y$ with $w(Y)=\tau$ for which the Boolean
algebras $RC(Y)$ and $RC(X)$ are isomorphic. Now, Proposition
\ref{absolute} implies that $X$ and $Y$ are co-absolute spaces.
\sqs

Obviously, if $X$ is co-absolute with a compact Hausdorff
 space $Y$ with $w(Y)=\tau$ then $\pi w(X)=\pi w(Y)\le\tau$.
Hence, we obtain:

\begin{cor}\label{pon2c1}{\rm (Ponomarev \cite{P1})}
A compact Hausdorff space $X$ is co-abso\-lu\-te with a compact
 metrizable space iff $\pi w(X)\le\aleph_0$.
\end{cor}

\section{On a problem of G. Birkhoff and some related problems. A
characterization of the spaces which are co-absolute with
(zero-dimensional) Eberlein compacts}

Recall that a space $X$ is called {\em semiregular} if $RO(X)$ is
a base for $X$.

\begin{notas}\label{bir}
\rm We will denote:
\begin{itemize}
\item by $\MM$ the class of all metrizable spaces,
\item by $\MM_0$ the class of all zero-dimensional metrizable spaces,
\item by $\MM_+$ the class of all regular Hausdorff (= $T_3$)
spaces $X$ which can be written in the form
$X=\bigoplus\{X_\g\st\g\in\GA\}$, where $\GA$ is an arbitrary set
and for every $\g\in\GA$, $w(X_\g)\le\aleph_0$,
\item by $\RR(\tau)$ the class of all $T_3$-spaces $X$ with
$w(X)=\tau$,
\item by $\SS\RR(\tau)$ the class of all semiregular spaces $X$ with
$w(X)=\tau$,
\item by $\DD$ the class of all discrete spaces,
\item by $\KK(\tau)$  (resp., by  $\KK_0(\tau)$) the class of
all compact Hausdorff (resp., and zero-dimensional) spaces $X$ with $w(X)\le\tau$,
\item by $\EE$ the class of all Eberlein compacts
(= weakly compact subsets of Banach spaces),
\item by $\EE_0$ the class of all zero-dimensional Eberlein
compacts,
\item by $\SS$ (respectively, by $\CC\SS$) the class of all  spaces
(respectively, all compact spaces) which have a dense Eberlein
subspace (where $``$Eberlein space" means $``$a subspace of an
Eberlein compact"),
\item by $\ZZ(\tau)$ (resp., $\ZZ\KK(\tau)$)
the class of all zero-dimensional
Hausdorff (resp., and compact) spaces $X$ with  $w(X)=\tau$.
\end{itemize}

If $\CC$ is a class of topological spaces, we will set
$\BB\CC=\{A\st A$ is a Boolean algebra and there exists $X\in\CC$
such that $A$ is isomorphic to the Boolean algebra $RO(X) \}$.
\end{notas}

The Problem 72 of G. Birkhoff \cite{Bi} is the following:
characterize internally the elements of the class $\BB\MM$. It was
solved by V. I. Ponomarev \cite{P1}. He proved the following
beautiful theorem: if $A$ is a complete Boolean algebra then
$A\in\BB\MM$ iff it has a $\s$-disjointed dense subset $B$ (i.e.
$B$ is a dense subset of $A$ and $B=\bigcup\{B_n\st
n\in\mathbb{N}^+\}$, where for every $n\in\mathbb{N}^+$ and for
every two different elements $a,b$ of $B_n$ we have $a\we b=0$).
The proof of this theorem is difficult. We will obtain a direct
(and easier) proof of it which leads to a characterization of the
class of spaces which are co-absolute with (zero-dimensional)
Eberlein compacts. Further, we will give some easily proved
solutions to some analogous problems. We will show that
$\BB\MM=\BB\EE$ and we will describe the elements of the classes
$\BB\MM_+$ and $\BB\ZZ(\tau)$ ($=\BB\ZZ\KK(\tau)=\BB\RR(\tau)$).
Clearly, $\DD\cup\RR(\aleph_0)\sbe\MM_+\sbe\MM$ and
$\KK(\aleph_0)\sbe\RR(\aleph_0)$. It is easy to see that the class
$\MM_+$ coincides with the class of all metrizable spaces which
have a metrizable locally compact extension. Note that if
$X\in\DD$ then $RO(X)=P(X)$; hence, by Tarski-Lindenbaum Theorem,
$A\in\BB\DD$ iff $A$ is a complete atomic Boolean algebra.

\begin{pro}\label{birk}
$\BB\MM=\BB\EE=\BB\CC\SS=\BB\SS$.
\end{pro}

\doc By a theorem of A. V. Arhangel'ski\u{i} \cite{Ar},
 every metric space can be
densely embedded in an Eberlein compact. Conversely, I. Namioka
\cite{Na} and Y. Benyamini-M. E. Rudin-M. Wage \cite{BRW} proved
that every Eberlein compact contains a dense  metrizable subspace.
Applying  \cite[Lemma 1.4]{D5}, we conclude that $\BB\MM=\BB\EE$.
Since every closed subset of an Eberlein compact is an Eberlein
compact, we get that $\BB\EE=\BB\CC\SS=\BB\SS$. \sqs

For proving the next theorem, we need to recall some facts and
definitions from \cite{D3,D4}.

\begin{defi}\label{alsb}{\rm \cite{D3,D4}}
\rm A family $\AA$ of subsets of a topological space $X$ is said
to be an {\em almost subbase of} $X$ if every element $V$ of $\AA$
has a representation $V=\bigcup\{U_n(V)\st n\in\mathbb{N}^+ \}$,
where for every $n\in\mathbb{N}^+ $, $U_n(V)\sbe U_{n+1}(V)$,
$U_{2n-1}(V)$ is a zero-set in $X$ and $U_{2n}(V)$ is a cozero-set
in $X$ (such a family $\{U_n(V)\st i\in\mathbb{N}^+\}$ will be
called an {\em Urysohn representation of} $V$), so that the family
$\AA\cup\{X\stm U_{2n-1}(V)\st V\in\AA, n\in\mathbb{N}^+ \}$ is a
subbase of $X$.
\end{defi}

\begin{theorem}\label{eberlein}{\rm \cite{D3,D4}}
A compact Hausdorff space is an Eberlein compact iff it has a
$\sigma$-point-finite almost subbase.
\end{theorem}

\begin{theorem}\label{birkhoffpon}
 A complete Boolean algebra $A$ is isomorphic to an algebra
of the form $RC(X)$, where $X$ is a (zero-dimensional) Eberlein
compact, iff $A$ has a $\s$-disjointed dense subset.
\end{theorem}

\doc ($\Rightarrow$) Let $A$ be a Boolean algebra which is
isomorphic to RC(X), where $X$ is an  Eberlein compact. As we have
already mentioned, there exists a metrizable dense subset $Y$ of
$X$. Hence $A$ is isomorphic to $RC(Y)$. The space $Y$ has a
$\sigma$-discrete base $\BB=\bigcup\{\BB_i\st i\in\mathbb{N}^+\}$,
where $\BB_i$ is a discrete family for every $i\in\mathbb{N}^+$.
Set, for every $i\in\mathbb{N}^+$, $\BB_i\ap=\{\cl(U)\st
U\in\BB_i\}$, and let $\BB\ap=\bigcup\{\BB_i\ap\st
i\in\mathbb{N}^+\}$. Then, obviously, $\BB\ap$ is a
$\sigma$-disjointed dense subset of $RC(Y)$. Hence, $A$  has a
$\s$-disjointed dense subset.

\noindent($\Leftarrow$) Let $A$ be a complete Boolean algebra
having a $\s$-disjointed dense subset $B_0$. Let $B$ be the
Boolean subalgebra of $A$ generated by $B_0$. Then $A$ is a
minimal completion of $B$. Set $X=S^a(B)$.
 Then $X$ is a
zero-dimensional compact Hausdorff space and there exists an
isomorphism $\p:B\lra CO(X)$. We will show that $\BB=\p(B_0)$ is a
$\s$-disjoint almost subbase of $X$.  For every $V\in\BB$ and
every $n\in\mathbb{N}^+$, set $U_n(V)=V$. Then $\{U_n(V)\st
i\in\mathbb{N}^+\}$ is an Urysohn representation of $V$. Hence, we
have to show that the family $\BB\ap=\BB\cup\{X\stm V\st
V\in\BB\}$ is a subbase of $X$. Obviously, $\BB\ap=\p(B_0\cup
B_0^*)$, where $B_0^*=\{b^*\st b\in B_0\}$. Since, clearly, the
set of all finite joins of all finite meets of the elements of the
subset $B_0\cup B_0^*$ of $A$ coincides with $B$, we get that the
family of all finite unions of the finite intersections of the
elements of the family $\BB\ap$ coincides with $CO(X)$ which is a
base of $X$. Hence, the family of all finite intersections of the
elements of $\BB\ap$ is a base of $X$, i.e. $\BB\ap$ is a subbase
of $X$. Therefore, $\BB$ is an almost subbase of $X$. Since $\BB$
is, obviously, a $\s$-disjoint family, we get, by Theorem
\ref{eberlein}, that $X$ is an Eberlein compact. Now, $RC(X)$ is a
minimal completion of $CO(X)$; thus  $RC(X)$ and $A$ are
isomorphic Boolean algebras. \sqs

Combining the last theorem with Proposition \ref{birk}, we obtain
the Ponomarev Theorem \cite{P1} giving a solution of Birkhoff's
Problem 72 \cite{Bi}.

\begin{cor}\label{birkhoffponom}{\rm (V. I. Ponomarev \cite{P1})}
 A complete Boolean algebra $A$ is isomorphic to an algebra
of the form $RC(X)$, where $X$ is a metrizable space, iff $A$ has
a $\s$-disjointed dense subset.
\end{cor}

Finally,  we get that:

\begin{cor}\label{birkhoffponomz}
$\BB\MM=\BB\CC\SS=\BB\SS=\BB\EE=\BB\EE_0=\BB\MM_0$.
\end{cor}

\doc By Proposition \ref{birk}, $\BB\MM=\BB\CC\SS=\BB\SS=\BB\EE$.
From Theorem \ref{birkhoffpon}, we get that $\BB\EE=\BB\EE_0$. Let
us prove  that $\BB\EE=\BB\MM_0$. Indeed, we have that
$\BB\MM_0\sbe\BB\MM=\BB\EE$.  Conversely, let $X$ be an Eberlein
compact. Then, by Theorem \ref{birkhoffpon}, there exists a
zero-dimensional Eberlein compact $Y$ such that $RC(X)\cong
RC(Y)$. Now, $Y$ has a dense metrizable subspace $Z$. Thus
$RC(X)\cong RC(Z)$ and $Z$ is a zero-dimensional metrizable space.
Therefore, $\BB\EE\sbe\BB\MM_0$. So, $\BB\EE=\BB\MM_0$. \sqs

\begin{theorem}\label{abseb}
Let $X$ be a compact Hausdorff space. Then the following
conditions are equivalent:

\noindent(a) $X$ is co-absolute with an Eberlein compact;

\noindent(b) $X$ has a $\s$-disjoint $\pi$-base;

\noindent(c)  $X$ is co-absolute with a zero-dimensional Eberlein
compact.
\end{theorem}

\doc (a)$\Rightarrow$(b) Let $Y$ be an Eberlein compact which is co-absolute
with $X$. Then $RC(Y)\cong RC(X)$. By Theorem \ref{birkhoffpon},
the Boolean algebra $RO(X)$ (which is isomorphic to the Boolean
algebra $RC(X)$) has a $\s$-disjointed dense subset $\AA$. Then,
obviously, $\AA$ is a $\s$-disjoint $\pi$-base of $X$.

(b)$\Rightarrow$(c) Let $\AA$ be a  $\s$-disjoint $\pi$-base of
$X$. Set $\AA\ap=\{\int(\cl(U))\st U\in\AA\}$. Then, obviously,
$\AA\ap$ is a $\s$-disjointed dense subset of the Boolean algebra
$RO(X)$. Since $RO(X)\cong RC(X)$, Theorem \ref{birkhoffpon}
implies that there exists a zero-dimensional Eberlein compact $Y$
with $RC(Y)\cong RC(X)$. Now Propositions \ref{absolute} and
\ref{absoluteun} imply that $X$ and $Y$ are co-absolute spaces.

(c)$\Rightarrow$(a) This is clear. \sqs

 We are now going to characterize  classes
 $\BB\ZZ(\tau)$($=$ $\BB\ZZ\KK(\tau)=\BB\SS\RR(\tau)$) and $\BB\MM_+$.

\begin{theorem}\label{birnul}
Let $A$ be a Boolean algebra and $\tau$ be an infinite cardinal
number. Then the following conditions are
equivalent:

\noindent(a) $A\in\BB\SS\RR(\tau)$;

\noindent(b) $A$ is complete and contains a dense subset $B$ with
$\card{B}=\tau$;

\noindent(c) $A\in\BB\ZZ\KK(\tau)$.

\noindent Note also that $\BB\ZZ(\tau)=\BB\ZZ\KK(\tau)$.
\end{theorem}

\doc (a)$\Rightarrow$(b) Let $A$ be isomorphic to $RO(X)$, where $X$ is a
semiregular space with $w(X)=\tau$.  There exists a subset $\BB$
of $RO(X)$ which is a base of $X$ and $\card{\BB}=\tau$. Then
$\BB$ is a dense subset of $RO(X)$.

\smallskip

\noindent(b)$\Rightarrow$(c) Let $A$ be a complete Boolean algebra
having a dense subset $B\ap$ with $\card{B\ap}=\tau$. Let $B$ be
the Boolean subalgebra of $A$ generated by $B\ap$.
 Then $\card{B}=\tau$ and $A$ is
a minimal completion of $B$. Set $X=S^a(B)$. Then $X$ is a compact
zero-dimensional Hausdorff space with $w(X)=\tau$. Since $B$ is
isomorphic to $CO(X)$ and $RC(X)$ is a minimal completion of
$CO(X)$, we get that $A$ is isomorphic to $RC(X)$. Hence,
$A\in\BB\ZZ\KK(\tau)$.

\smallskip

\noindent(c)$\Rightarrow$(a) This is obvious.
 \sqs

Note that the last assertion, the Brouwer topological
characterization of the Cantor set $C$ as the unique (up to
homeomorphism) dense in itself zero-dimensional compact metrizable
space and the obvious fact that the atoms of a Boolean algebra $A$
correspond to the isolated points of the dual spaces of the LCAs
of the form $(A,\rho,\BBBB)$ imply the second algebraic
characterization of $RC(C)$ mentioned above  (namely, that $RC(C)$
is the unique (up to isomorphism) atomless complete Boolean
algebra containing a countable dense subalgebra).

\begin{lm}\label{birkhoff}
 If $A$ is a Boolean algebra then $A\in\BB\MM_+$   if and only if
   $A=\prod\{A_\g\st \g\in\GA\}$
where, for every $\g\in\GA$, $A_\g$ is a complete Boolean algebra
and there exists a normal contact relation $C_\g$ on $A_\g$ such
that $w(A,C_\g)\le\aleph_0$.
\end{lm}

\doc ($\Rightarrow$) Let $A$ be isomorphic to $RO(X)$ for some
$X\in\MM_+$.  Since the Boolean algebras $RO(X)$ and $RC(X)$ are
isomorphic, we get that $A$ is isomorphic to $RC(X)$. We have that
$X=\bigoplus\{X_\g\st\g\in\GA\}$, where $\GA$ is a set and for
every $\g\in\GA$, $w(X_\g)\le\aleph_0$.  The spaces $X_\g$,
$\g\in\GA$, are metrizable; hence they have  metrizable
compactifications $cX_\g$. Then $L=\bigoplus\{cX_\g\st \g\in\GA\}$
is a (metrizable) locally compact extension of $X$ and, by
\cite[Lemma 1.4]{D5}, $RC(X)$ is isomorphic to $RC(L)$. So, by
\cite[Theorem 2.14]{D5}, \ref{prodclcacat} and \ref{weightpiw},
$RC(X)$ is isomorphic to $\prod\{RC(cX_\g)\st \g\in\GA\}$, where
$w(RC(cX_\g),\rho_{cX_\g})=w(cX_\g)\le\aleph_0$  (see
\cite[Example 1.3]{D5} for $\rho_{cX_\g}$).

\smallskip

\noindent($\Leftarrow$) Let $A=\prod\{A_\g\st \g\in\GA\}$ where,
for every $\g\in\GA$, $A_\g$ is a complete Boolean algebra and
there exists a normal contact relation $C_\g$ on $A_\g$ such that
$w(A,C_\g)\le\aleph_0$. Let $(A,\rho,\BBBB)=\prod\{(A_\g,C_\g)\st
\g\in\GA\}$ (see Definition \ref{prodclca}). Then, by
\ref{prodclcacat}, $\{\pi_\g:(A,\rho,\BBBB)\lra
(A_\g,C_\g)\st\g\in\GA\}$ is a $\DLC$-product of the family
$\{(A_\g,C_\g)\st \g\in\GA\}$. Let, for every $\g\in\GA$,
$X_\g=\LAM^a(A_\g,C_\g)$ and $X=\LAM^a(A,\rho,\BBBB)$. Then, by
\ref{prodclcacat} and \cite[Theorem 2.14]{D5},
$X=\bigoplus\{X_\g\st\g\in\GA\}$ and, by \ref{weightpiw},
$w(X_\g)\le\aleph_0$ for every $\g\in\GA$.   Hence $X\in\MM_+$
and, by  \cite[Theorem 2.14]{D5}, $A$ is isomorphic to $RO(X)$.
\sqs

By  \cite[Example 1.2]{D5}, for every complete Boolean algebra
$B$, $(B,\rho_s)$ is a CNCA. If $w(B,\rho_s)\le\aleph_0$ then
$\card{B}=w(B,\rho_s)\le\aleph_0$ (because, by \cite[Example
1.2]{D5}, for every $b\in B$, we have that $b\ll_{\rho_s} b$).
Since $B$ is complete, it follows that $B$ is a finite Boolean
algebra (see, e.g., \cite{Si}), and hence $B=\2^n$ for some
$n\in\mathbb{N}^+$. Therefore, if in Lemma \ref{birkhoff} we set,
for every $\g\in\GA$, $C_\g=\rho_s$ then we will obtain that
$A=\2^{\card{\GA}}$, i.e. that $A$ is a complete atomic Boolean
algebra.

\begin{theorem}\label{m0c}
Let $A$ be a Boolean algebra. Then $A\in\BB\MM_+$   iff
$A=\prod\{A_\g\st \g\in\GA\}$ where, for every $\g\in\GA$, $A_\g$
is a complete Boolean algebra having a dense countable subset.
\end{theorem}

\doc It follows from  \ref{birkhoff} and \ref{nuldimlm}. \sqs

\section{A completion theorem for LCAs}

\begin{defi}\label{defcompl}
\rm Let $(A,\rho,\BBBB)$ be an LCA. A pair
$(\p,(A\ap,\rho\ap,\BBBB\ap))$ is called an {\em LCA-completion}
of the LCA $(A,\rho,\BBBB)$ if $(A\ap,\rho\ap,\BBBB\ap)$ is a
CLCA, $\p$ is an LCA-embedding (see \cite[Definition 1.11]{D5}) of
$(A,\rho,\BBBB)$ in $(A\ap,\rho\ap,\BBBB\ap)$, and $\p(\BBBB)$ is
a dV-dense subset of $(A\ap,\rho\ap,\BBBB\ap)$ (see \ref{clcawe}
for the last notion).

Two LCA-completions $(\p,(A\ap,\rho\ap,\BBBB\ap))$ and
$(\psi,(A'',\rho'',\BBBB''))$ of a local contact algebra
$(A,\rho,\BBBB)$ are said to be {\em equivalent} if there exists
an LCA-isomorphism
$\eta:(A\ap,\rho\ap,\BBBB\ap)\lra(A'',\rho'',\BBBB'')$ such that
$\psi=\eta\circ\p$.
\end{defi}

Note that condition (BC3) (see  \cite[Definition 1.11]{D5})
implies that every dV-dense subset of an LCA $(A,\rho,\BBBB)$ is a
dense subset of $A$. Hence, if $(\p,(A\ap,\rho\ap,\BBBB\ap))$ is
an LCA-completion of the LCA $(A,\rho,\BBBB)$ then $(\p,A\ap)$ is
a minimal completion of the Boolean algebra $A$.

 Let us start with a simple lemma.

 \begin{lm}\label{lmcompletion}
Let $(\p,(B,\eta,\BBBB\ap))$ be an LCA-completion of an LCA
$(A,\rho,\BBBB)$ and let us suppose, for simplicity, that $A\sbe
B$ and $\p(a)=a$ for every $a\in A$. Then:

\smallskip

\noindent(a) $\BBBB\ap=\downarrow_B(\BBBB)$ and $\BBBB\ap\cap
A=\BBBB$;

\smallskip

\noindent(b) If $J$ is a $\d$-ideal (see \cite[Definition
2.2]{D5}) of $(B,\eta,\BBBB\ap)$ then $J\cap A$ is a $\d$-ideal of
$(A,\rho,\BBBB)$ and $\downarrow_B(J\cap A)=J$;

\smallskip

\noindent(c) If $J$ is a $\d$-ideal  of $(A,\rho,\BBBB)$ then
$\downarrow_B(J)$ is a $\d$-ideal of $(B,\eta,\BBBB\ap)$ and
$A\cap\downarrow_B(J)=J$;

\smallskip

\noindent(d) If $J$ is a prime element (see the text before
\cite[Proposition 3.4]{D5}) of $I(B,\eta,\BBBB\ap)$ (see
\cite[Definition 2.2]{D5}) then $J\cap A$ is a prime element of
the frame $I(A,\rho,\BBBB)$;

\noindent(e) If $J$ is a prime element  of $I(A,\rho,\BBBB)$ then
$\downarrow_B(J)$ is a prime element of $(B,\eta,\BBBB\ap)$.
\end{lm}

\doc (a)   Let
$b\in\BBBB\ap$. Then, by condition (BC1) (see \cite[Definition
1.11]{D5}), there exists $c\in\BBBB\ap$ such that $b\lle c$
(because $b\lle 1$). Since $\BBBB$ is a dV-dense subset of
$(B,\eta,\BBBB\ap)$, there exists $a\in\BBBB$ such that $b\le a\le
c$.  Hence $\BBBB\ap\sbe\downarrow_B(\BBBB)$. Since
$\BBBB\sbe\BBBB\ap$ and $\BBBB\ap$ is an ideal of $B$, we get that
$\downarrow_B(\BBBB)\sbe\BBBB\ap$. Hence,
$\BBBB\ap=\downarrow_B(\BBBB)$.

Obviously, $\BBBB\sbe\BBBB\ap\cap A$. If $a\in\BBBB\ap\cap A$
then, as above, there exists $b\in\BBBB$ such that $a\le b$. Thus
$a\in\BBBB$. Hence $\BBBB\ap\cap A=\BBBB$.

\smallskip

\noindent(b) We have that $J\cap A\sbe \BBBB\ap\cap A=\BBBB$.
 Let $a\in J\cap A$. Then there
exists $b\in J$ such that $a\lle b$. Since $\BBBB$ is a dV-dense
subset of $(B,\eta,\BBBB\ap)$, we get that there exists
$c\in\BBBB$ such that $a\lle c\lle b$ (see Fact \ref{wfact}). Then
$c\in J\cap A$ and $a\llx c$. So, $J\cap A$ is a $\d$-ideal of
$(A,\rho,\BBBB)$. The last argument shows as well that
$J\sbe\downarrow_B(J\cap A)$. Since, clearly, $\downarrow_B(J\cap
A)\sbe J$, we get that $\downarrow_B(J\cap A)= J$.

\smallskip

\noindent(c) Let  $J$ be a $\d$-ideal  of $(A,\rho,\BBBB)$. Set
$J\ap=\downarrow_B(J)$. Clearly, $J\ap$ is an ideal of $B$. Let
$a\in J\ap$. Then there exists $b,c\in J$ such that $a\le b\llx
c$. Thus $a\lle c$ and $c\in J\ap$. Hence $J\ap$ is a $\d$-ideal
of $(B,\eta,\BBBB\ap)$. Obviously, $J\sbe A\cap\downarrow_B(J)$.
Conversely, let $a\in  A\cap\downarrow_B(J)$. Then there exists
$b\in J$ such that $a\le b$. Thus $a\in J$. So,
$A\cap\downarrow_B(J)=J$.

\smallskip

\noindent(d) Let $J$ be a prime element  of $I(B,\eta,\BBBB\ap)$.
 Then, by (b), $J\cap A\in I(A,\rho,\BBBB)$. Let $J_1,J_2\in
 I(A,\rho,\BBBB)$ and $J_1\cap J_2\sbe J\cap A$. Then
 $\downarrow_B(J_1)\cap\downarrow_B(J_2)=$ $\downarrow_B(J_1\cap
 J_2)\sbe\downarrow_B(J\cap A)$. Since, by (c), $\downarrow_B(J_i)\in
 I(B,\eta,\BBBB\ap)$, for $i=1,2$, and, by (b), $\downarrow_B(J\cap A)=J$,
  we get that $\downarrow_B(J_1)\sbe
 J$ or $\downarrow_B(J_2)\sbe
 J$. Then $A\cap\downarrow_B(J_1)\sbe
 A\cap J$ or $A\cap\downarrow_B(J_2)\sbe
 A\cap J$. Thus, by (c), $J_1\sbe J\cap A$ or $J_2\sbe J\cap A$.
 Hence, $J\cap A$ is a prime element of $I(A,\rho,\BBBB)$.

\smallskip

\noindent(e) Let $J$ be a prime element of $I(A,\rho,\BBBB)$. Let
$J_1,J_2\in I(B,\eta,\BBBB\ap)$ and  $J_1\cap J_2\sbe
\downarrow_B(J)$. Then, by (c), $A\cap J_1\cap J_2\sbe
A\cap\downarrow_B(J)=J$. Hence, by (b), $A\cap J_1\sbe J$ or
$A\cap J_2\sbe J$. Thus, by (b), $J_1\sbe\downarrow_B(J)$ or
$J_2\sbe\downarrow_B(J)$. Therefore, $\downarrow_B(J)$ is a prime
element of $I(B,\eta,\BBBB\ap)$. \sqs

\begin{theorem}\label{lcacompletion}
Every LCA $(A,\rho,\BBBB)$ has a unique (up to equivalence)
LCA-completion.
\end{theorem}

\doc Let  $(A,\rho,\BBBB)$ be an LCA. Then, by Roeper's theorem
\cite[Theorem 2.1]{D5}, there exists a locally compact Hausdorff
space $X$ and an LCA embedding
$\lag:(A,\rho,\BBBB)\lra(RC(X),\rho_X,CR(X))$ such that
$\{\int(\lag(a))\st a\in\BBBB\}$ is a base of $X$. Since $\BBBB$
is closed under finite joins, we get easily (using the compactness
of the elements of $CR(X)$) that $\lag(\BBBB)$ is a dV-dense
subset of the CLCA $(RC(X),\rho_X,CR(X))$. Hence the pair $$(\lag,
(RC(X),\rho_X,CR(X)))$$ is an LCA-completion of the LCA
$(A,\rho,\BBBB)$.

We will now prove the uniqueness (up to equivalence) of the
LCA-completion. Let  $(\p,(B,\eta,\BBBB\ap))$ be an LCA-completion
of the LCA $(A,\rho,\BBBB)$. Then, as we have already mentioned,
$(\p,B)$ is a minimal completion of $A$, i.e. the Boolean algebra
$B$ is determined uniquely (up to isomorphism)  by the Boolean
algebra $A$. We can suppose wlog that $A\sbe B$ and $\p(a)=a$, for
every $a\in A$. Thus $A$ is a Boolean subalgebra of $B$.

As we have already shown (see Lemma \ref{lmcompletion}(a)),
$\BBBB\ap=\downarrow_B(\BBBB)$, i.e. the set $\BBBB\ap$ is
uniquely determined by the set $\BBBB$.

We have that $\eta_{|A}=\rho$. We will show that the relation
$\eta$ on $B$ is uniquely determined by the relation $\rho$ on
$A$. There are two cases.

\smallskip

\noindent{\em Case 1.} Let $a_1\in\BBBB\ap$ and $b_1\in B$. We
will prove that $a_1\lle b_1$ iff there exist $a,b\in\BBBB$ such
that $a_1\le a\llx b\le b_1$. By (BC1), it is enough to prove this
for $b_1\in\BBBB\ap$.

So, let
 $a_1,b_1\in\BBBB\ap$ and $a_1\lle b_1$. Then,
 using dV-density of $\BBBB$ in $(B,\eta,\BBBB\ap)$ and Fact \ref{wfact},
  we get that there exist $a,b\in\BBBB$ such that
 $a_1\le a\lle b\le b_1$. Then $a\llx b$.

   The converse assertion is clear because, for
every $a,b\in A$, $a\llx b$ iff $a\lle b$.

\smallskip

\noindent{\em Case 2.}  Let $a_1\in B\stm\BBBB\ap$ and $b_1\in B$.
We will prove that $a_1\lle b_1$ iff (for every prime element $J$
of $I(A,\rho,\BBBB)$) [(there exists $a\in \downarrow_B(\BBBB)\stm
\downarrow_B(J)$ such that $a\lle a_1^*$) or (there exists $b\in
\downarrow_B(\BBBB)\stm \downarrow_B(J)$ such that $b\lle b_1$)].
Note that the inequalities $a\lle a_1^*$ and $b\lle b_1$ from the
above formula  are already expressed in Case 1 in a form which
depends only of $(A,\rho,\BBBB)$ (because $a,b\in\BBBB\ap$).
Hence,
  Case 1 and Case 2 will imply that the relation $\eta$ on
$B$ is uniquely determined by the relation $\rho$ on $A$.

So,  let $a_1\in B\stm\BBBB\ap$ and $b_1\in B$. Then using
\cite[(25)]{D5}, \cite[Proposition 3.4]{D5}, \cite[Proposition
3.6]{D5}, and Lemma \ref{lmcompletion}, we get that $a_1\lle b_1$
iff $a_1(-\eta) b_1^*$ iff [(for every
$\s\in\Psi^a(B,\eta,\BBBB\ap)$)($\{a_1,b_1^*\}\not\sbe\s$)] iff
(for every prime element $J\ap$ of $I(B,\eta,\BBBB\ap)$)[(there
exists $a\in\BBBB\ap\stm J\ap$ such that $a(-\eta)a_1$) or (there
exists $b\in\BBBB\ap\stm J\ap$ such that $b(-\eta)b_1^*$)]
 iff (for every prime element $J$ of
$I(A,\rho,\BBBB)$) [(there exists $a\in \downarrow_B(\BBBB)\stm
\downarrow_B(J)$ such that $a\lle a_1^*$) or (there exists $b\in$
 $\downarrow_B(\BBBB)\stm \downarrow_B(J)$ such that $b\lle b_1$)].

Let now $(\p_1,(A_1,\rho_1,\BBBB_1))$ and
$(\p_2,(A_2,\rho_2,\BBBB_2))$ be two LCA-comple\-tions of an LCA
$(A,\rho,\BBBB)$. Then, since $(\p_i,A_i)$, for $i=1,2$, are
minimal completions of $A$, there exists a Boolean isomorphism
$\p:A_1\lra A_2$ such that $\p\circ\p_1=\p_2$. The preceding
considerations imply that $\BBBB_i=\downarrow_{A_i}(\p_i(\BBBB))$,
for $i=1,2$. From this we easily get that $\p(\BBBB_1)=\BBBB_2$.
Further, for $a_i\in\BBBB_i,b_i\in A_i$, $i=1,2$, we have that
$a_i\ll_{\rho_i} b_i$ iff there exists $a_i\ap,b_i\ap\in\BBBB$
such that $a_i\ap\llx b_i\ap$, $a_i\le\p_i(a_i\ap)$ and
$\p_i(b_i\ap)\le b_i$, for $i=1,2$. Finally, for $a_i\in
A_i\stm\BBBB_i,b_i\in A_i$, $i=1,2$, we have that $a_i\ll_{\rho_i}
b_i$ iff  (for every prime element $J$ of $I(A,\rho,\BBBB)$)
[(there exists $a_i\ap\in \BBBB_i\stm \downarrow_{A_i}(\p_i(J))$
such that $a_i\ap\ll_{\rho_i} a_i^*$) or (there exists $b_i\ap\in
\BBBB_i\stm \downarrow_{A_i}(\p_i(J))$ such that
$b_i\ap\ll_{\rho_i} b_i$)]. Having in mind these formulas, it is
easy to conclude that $\p$ is an LCA-isomorphism. Hence the
LCA-completions $(\p_1,(A_1,\rho_1,\BBBB_1))$ and
$(\p_2,(A_2,\rho_2,\BBBB_2))$ of $(A,\rho,\BBBB)$ are equivalent.
\sqs

\begin{cor}\label{corcompl}
Let $(A,\rho,\BBBB)$ be an LCA and $(B,\eta,\BBBB\ap)$ be a CLCA.
Then $\Psi^a(A,\rho,\BBBB)$ is homeomorphic to
$\Psi^a(B,\eta,\BBBB\ap)$ if and only if there exists an
LCA-embedding $\p:(A,\rho,\BBBB)\lra(B,\eta,\BBBB\ap)$ such that
$\p(\BBBB)$ is a dV-dense subset of $(B,\eta,\BBBB\ap)$.
\end{cor}

\doc ($\Rightarrow$) In the proof of Theorem \ref{lcacompletion},
we have seen that the set $\lag(\BBBB)$ is dV-dense in
$\Psi^t(\Psi^a(A,\rho,\BBBB))$. Since
$\Psi^t(\Psi^a(A,\rho,\BBBB))$ is LCA-isomorphic to
$\Psi^t(\Psi^a(B,\eta,\BBBB\ap))$ and
$(B,\eta,\BBBB\ap)\cong\Psi^t(\Psi^a(B,\eta,\BBBB\ap))$, we get
that there exists an LCA-embedding
$\p:(A,\rho,\BBBB)\lra(B,\eta,\BBBB\ap)$ such that $\p(\BBBB)$ is
a dV-dense subset of $(B,\eta,\BBBB\ap)$.

\noindent($\Leftarrow$) By the proof of Theorem
\ref{lcacompletion}, $(\lag,\Psi^t(\Psi^a(A,\rho,\BBBB)))$ is an
LCA-comple\-tion of $(A,\rho,\BBBB)$. Since the hypothesis of our
assertion imply that the pair $(\p,(B,\eta,\BBBB\ap))$ is also an
LCA-completion of $(A,\rho,\BBBB)$, we get, by Theorem
\ref{lcacompletion}, that the CLCAs $\Psi^t(\Psi^a(A,\rho,\BBBB))$
and $(B,\eta,\BBBB\ap)$ are LCA-isomorphic. Then
$\Psi^a(B,\eta,\BBBB\ap)\cong\Psi^a(\Psi^t(\Psi^a(A,\rho,\BBBB)))\cong\Psi^a(A,\rho,\BBBB)$.
\sqs

\begin{cor}\label{corcompl1}
Let $(A,\rho,\BBBB)$ and $(A\ap,\rho\ap,\BBBB\ap)$ be LCAs. Then
$\Psi^a(A,\rho,\BBBB)$ is homeomorphic to
$\Psi^a(A\ap,\rho\ap,\BBBB\ap)$ iff there exists a CLCA
$(B,\eta,\BBBB'')$ and  LCA-embeddings
$\p:(A,\rho,\BBBB)\lra(B,\eta,\BBBB'')$ and
$\p\ap:(A\ap,\rho\ap,\BBBB\ap)\lra(B,\eta,\BBBB'')$ such that the
sets $\p(\BBBB)$ and $\p\ap(\BBBB\ap)$ are dV-dense in
$(B,\eta,\BBBB'')$.
\end{cor}

\doc ($\Rightarrow$) Set $(B,\eta,\BBBB'')=\Psi^t(\Psi^a(A,\rho,\BBBB))$. Then, by  the hypothesis of our
assertion, there exists an LCA-isomorphism
$\psi:\Psi^t(\Psi^a(A\ap,\rho\ap,\BBBB\ap))\lra(B,\eta,\BBBB'')$.
Now, it is clear that the maps $\lag$ and $\psi\circ\l_{A\ap}^g$
are the required LCA-embeddings.

\noindent($\Leftarrow$) By Corollary \ref{corcompl}, we have that
$$\Psi^a(A,\rho,\BBBB)\cong\Psi^a(B,\eta,\BBBB'')\cong\Psi^a(A\ap,\rho\ap,\BBBB\ap).$$
\sqs

\baselineskip = 1.2\normalbaselineskip

\baselineskip = 1.00\normalbaselineskip

\vspace{0.25cm}

 Faculty  of Mathematics and Informatics,

Sofia University,

5 J. Bourchier Blvd.,

1164 Sofia,

Bulgaria


\begin{thebibliography}{99}
{\small

\bibitem{AHS}
  J. Ad\'amek, H. Herrlich, G. E. Strecker,
  \newblock  {\em Abstract and Concrete Categories},
   \newblock Wiley Interscience, New York, 1990.


\bibitem{A}
  P. S. Alexandroff,
  \newblock  {\em Outline of  set theory and general topology},
   \newblock Nauka, Moskow, 1977 (in Russian).


\bibitem{Ar}
A. V. Arhangel'ski\v{i},
\newblock {\em Compact Hausdorff spaces and unions of countable families of
metrizable spaces},
  \newblock Dokl. Akad. Nauk S.S.S.R., 232 (1977), 989--992 (in Russian);
   Soviet Math. Dokl.,  18 (1977), 165--169.


\bibitem{BRW}
 Y. Benyamini, M. E. Rudin, M. L. Wage,
\newblock {\em Continuous images of weakly compact subsets of Banach spaces},
  \newblock  Pacific J. Math.,  70(2) (1977), 309--324.




\bibitem{Bi}
  G. Birkhoff,
  \newblock  {\em Lattice Theory},
   \newblock Providence, Rhode Island, 1967.




\bibitem{dL}
T. de Laguna,
 \newblock {\em Point, line and surface as sets of solids},
 \newblock The Journal of Philosophy, 19 (1922), 449--461.

\bibitem{dV}
H. de Vries,
\newblock  {\em Compact Spaces and Compactifications, an Algebraic Approach},
\newblock Van Gorcum, The Netherlands, 1962.


\bibitem{D5}
G. Dimov,
\newblock  {\em A  De Vries-type Duality Theorem for locally compact spaces -- I},
\newblock arXiv:0903.2589v2, 1-37.



\bibitem{D}
G. Dimov,
\newblock  {\em A generalization of De Vries Duality Theorem},
\newblock Applied Categorical Structures (to appear) (DOI 10.1007/s10485-008-9144-5).

\bibitem{D1}
G. Dimov,
\newblock  {\em Some generalizations of Fedorchuk Duality Theorem -- I},
\newblock Topo\-logy Appl. 156 (2009), 728-746.

\bibitem{D2}
G. Dimov,
\newblock  {\em Some generalizations of Fedorchuk Duality Theorem -- II},
\newblock arXiv:0710.0181v1, 1-20.


\bibitem{D3}
G. Dimov,
\newblock  {\em On Eberlein spaces and related spaces},
\newblock C. R. Acad. Sc. Paris 304 (9) (S\'{e}rie I) (1987), 233--235.

\bibitem{D4}
G. Dimov,
\newblock  {\em Baire subspaces of} $c_0(\Gamma)$ {\em have dense}
$G_\d$ {\em metrizable subsets},
\newblock Rendiconti
Circolo Matemat. Palermo, Ser. II, Suppl. no.18 (1988), 275--285.







\bibitem{Dw}
 Ph. Dwinger,
\newblock  {\em Introduction to Boolean Algebras},
\newblock Physica Verlag, W\"{u}rzburg, 1961.





\bibitem{E}
 R. Engelking,
\newblock  {\em General Topology},
\newblock PWN, Warszawa, 1977.


\bibitem{F}
 V. V. Fedorchuk,
\newblock  {\em Boolean $\d$-algebras and quasi-open mappings},
\newblock  Sibirsk. Mat. \v{Z}. 14 (5) (1973), 1088--1099; English translation: Siberian Math. J.
14 (1973), 759-767 (1974).



\bibitem{J}
P. T. Johnstone,
\newblock {\em Stone Spaces},
\newblock Cambridge Univ. Press, Cambridge, 1982.




\bibitem{LE}
  S. Leader,
\newblock {\em Local proximity spaces},
\newblock  Math. Annalen 169 (1967), 275--281.



\bibitem{MR}
J. Mioduszewski, L. Rudolf,
\newblock H-closed and extremally disconected Hausdorff spaces,
\newblock  Dissert. Math. (Rozpr. Mat.) 66 (1969) 1--52.



\bibitem{Na}
I. Namioka,
\newblock {\em Separate continuity and joint continuity},
  \newblock  Pacific J. Math.,  51 (1974), 515--531.


\bibitem{P}
V. I. Ponomarev,
\newblock Paracompacta: their projection spectra and continuous
mappings,
\newblock   Mat. Sb. (N.S.) 60 (102) (1963) 89--119. (In Russian)

\bibitem{P1}
V. I. Ponomarev,
\newblock  Spaces  co-absolute with metric spaces,
\newblock   Uspekhi Mat. Nauk. 21 (4) (1966) 101--132; English translation: Russian
Math. Surveys 21 (4) (1966), 87-114.

\bibitem{PS}
V. I. Ponomarev, L. B. \v{S}apiro,
\newblock Absolutes of topological spaces and their continuous
maps,
\newblock Uspekhi Mat. Nauk 31 (1976) 121--136; English translation: Russian
Math. Surveys 21 (1976) 138-154.



\bibitem{R}
P. Roeper,
\newblock  {\em Region-based topology},
\newblock  Journal of Philosophical Logic 26 (1997), 251--309.

\bibitem{Si}
  R. Sikorski,
\newblock {\em Boolean Algebras},
\newblock Springer-Verlag, Berlin, 1964.

\bibitem{SS}
 L. A. Steen, J. A. Seebach,
\newblock  {\em Counterexamples in Topology},
\newblock Springer-Verlag, Berlin, 1978.


\bibitem{ST1}
 M. H. Stone,
\newblock  {\em  Theory of representation for Boolean algebras},
\newblock  Trans. Amer. Math. Soc. 40 (1936), 37-111.





\bibitem{ST}
 M. H. Stone,
\newblock  {\em Applications of the theory of Boolean rings to general
topology},
\newblock  Trans. Amer. Math. Soc. 41 (1937), 375-481.
}
\end{thebibliography}
\end{document}